\documentclass[review]{elsarticle}

\usepackage{lineno,hyperref}

\usepackage{graphicx}
\usepackage{epstopdf}
\usepackage{amsmath,amsfonts,amsthm,bm} 
\usepackage{soul,color} 

\definecolor{orange}{rgb}{1,0.5,0}
\newcommand{\orange}[1]{\textcolor{orange}{#1}}  
\newcommand{\blue}[1]{\textcolor{blue}{#1}}
\newcommand{\green}[1]{\textcolor{green}{#1}}
\newcommand{\red}[1]{\textcolor{red}{#1}}
\usepackage{multirow}
\DeclareGraphicsExtensions{.eps,.pdf,.png,.jpg,.jpeg,.xcf}
\usepackage{subcaption}
\usepackage[justification=centering]{caption}
\usepackage{placeins}
\usepackage{algorithm}
\usepackage{algpseudocode}
\modulolinenumbers[5]

\journal{Journal of \LaTeX\ Templates}

\begin{document}

\begin{frontmatter}

\title{A projection-based numerical integration scheme for embedded interface: Application to fluid-structure interaction}

\author[mymainaddress]{B. Liu\corref{mycorrespondingauthor}}
\ead{a0098961@u.nus.edu}

\author[mysecondaryaddress]{R. K. Jaiman}

\author[mymainaddress]{D. Tan \corref{mycorrespondingauthor}}
\cortext[mycorrespondingauthor]{Corresponding author}
\ead{mpetds@nus.edu.sg}

\address[mymainaddress]{Department of Mechanical Engineering, National University of Singapore, Singapore}
\address[mysecondaryaddress]{Department of Mechanical Engineering, University of British Columbia, Vancouver, Canada}

\begin{abstract} 

We present a projection-based numerical integration technique to deal with embedded interface in finite element (FE) framework. The element cut by an embedded interface is denoted as a cut cell. We recognize elemental matrices of a cut cell can be reconstructed from the elemental matrices of its sub-divided cells, via projection at matrix level. These sub-divided cells are termed as integration cells. The proposed technique possesses following characteristics (1) no change in FE formulation and quadrature rule; (2) consistency with the derivation of FE formulation in variational principle. It can be considered as a re-projection of the residuals of equation system in the test function space or a reduced-order modeling (ROM) technique.
These characteristics significantly improves its scalability, easy-to-implementation and robustness to deal with problems involving embedded discontinuities in FE framework. Numerical examples, e.g., vortex-induced vibration (VIV), rotation, free fall and rigid-body contact in which the proposed technique is implemented to integrate the variational form of Navier-Stokes equations in cut cells, are presented. 
\end{abstract}


\end{frontmatter}

\section{Introduction}

The embedded interface FE formulation is an appealing approach in problems involving moving interfaces, e.g., fluid-structure interaction or free surface flows, or situations in which efforts are made to eliminate the generation of body-fitted meshes. A number of schemes, overall termed as unfitted finite element methods, were proposed to weakly imposed boundary conditions along the embedded interfaces. For instance, partition of unity method~\cite{melenk1996partition,Babuska1997Ijfnmie}, eXtended/generalized finite element approach~\cite{Moees1999Ijfnmie,Duarte2000CS,Strouboulis2000Cmiamae,Strouboulis2001Cmiamae,Sukumar2001Cmiamae,Belytschko1999Ijfnmie}, fictitious domain method (FDM)~\cite{Ramiere2007JoCP,Glowinski2007CMiAMaE} and FCM, which is a combination of the fictitious domain technique and the high-order finite element approach~\cite{parvizian2007finite,Duester2008Cmiamae}. 
 
The concept of distributing Lagrange points along the embedded interfaces was well established and implemented in the aforementioned numerical approaches. However, an appropriate choice of Lagrange multiplier basis space is critical to satisfy the Bab{\v u}ska-Brezzi (BB) condition \cite{brezzi1974existence,brezzi1990discourse}. Recently, Nitsche's method~\cite{nitsche1971variationsprinzip} gained attention among research community due to its advantageous characteristics, e.g., variationally consistency and no increment in system size. It had been implemented to investigate a number of fluid-structure-interaction (FSI) problems, e.g., \cite{Burman2012ANM,Massing2014JoSC,dettmer2016stabilised,schillinger2016non,kadapa2017stabilised,zou2017nitsche}. In the present work, we implement Nitsche's method to weakly impose the Dirichlet boundary along the embedded interface. 
  
In all unfitted interface formulations, the numerical integration over the cut cell requires a special attention. Without appropriate treatments to ensure accurate approximations around the interface, the idea of unfitted finite element method becomes impractical. Overall, five important classes of integration methods consisting of embedded discontinuity in finite element formulation can be listed as (1) tessellation, (2) moment fitting methods, (3) methods based on the divergence theorem, (4) equivalent polynomial and (5) conformal mapping. Tessellation~\cite{liu2002mesh,belytschko2009review} is a well-established method, in which the cut cell is triangulated or quadrangulated into smaller integration cells. Its advantage is the embedded discontinuity can be accurately captured by aligning with the edge of integration cells. On the other hand, an uniform refinement~\cite{parvizian2007finite} of the cut cell can be implemented, to avoid the difficulty in aligning the embedded discontinuities with integration cells for complex geometries. 
However, the uniform refinement approach is computational expensive to obtain sufficient accurate numerical results. To improve the computational efficiency, adaptive refinement techniques, e.g., Quatree or smart Octree~\cite{samet1990applications,de2000computational}, can be used to minimize the integration error around the embedded interface. Nonetheless, the computational cost is still high compared with tessellation. Another improvement is to modify the integration weights of the standard Gauss quadrature \cite{rabczuk2007meshfree}, in which the weights are scaled based on the ratio of cut area by the discontinuity. The recent development in numerical integration techniques focuses on the elimination of subdivision on the cut cell, e.g., finding equivalent polynomial functions~\cite{ventura2006elimination}, constructing efficient quadrature rules for individual integration cell (moment-fitting
equations)~\cite{muller2013highly}, transforming the volume/surface integral to surface/line integral (divergence theorem)~\cite{hubrich2015numerical} and Schwarz-Christoffel conformal mapping~\cite{Natarajan2010IJfNMiE}. 

In the present work, we propose a projection-based numerical integration technique for problems in unfitted FE formulation. It works for both tessellation and adaptive refinement methods. Of particular we address the following issues: (1) easy-of-implementation in FE formulations (simplicity and scalibility), (2) capable to produce accurate numerical results (accuracy), (3) well-suited for FE formulation (variationally consistency) and (4) applicable to various FSI applications (robustness). The proposed numerical integration technique is a variant of tessellation method. In terms of algorithm, the primary differences are (1) no change in FE formulation and quadrature rule for the elements with/without embedded discontinuities, (2) the elemental matrices of subdivided integration cell are assembled via transformation in a quadratic form, a projection procedure. The transformation operation refers to the operation of changing the representation of a matrix between different bases a transformation tensor, such that the matrix retains equivalent. In finite element theorem, this assembly procedure is rooted in Bubnov-Galerkin method, a re-projection of residuals of equation system in the test function space. Alternatively, it can be considered as a reduced-order modeling (ROM) technique, where a higher dimension problem is projected into a lower dimension space. It is proven in present work the reconstructed elemental matrices via our proposed technique exactly recover the original elemental matrices obtained by the standard Gauss quadrature on the cut cell. 

The manuscript is organized as follows. The mathematical formulation of proposed numerical integration technique is discussed at first in Sect.~\ref{sec:NI}. The governing equations and FSI schemes are listed in Sect.~\ref{sec:govern}. The complete variational formulation of our unfitted FSI solver is shown in Sect.~\ref{sec:formulation}. Following that, the error analysis of this proposed numerical integration technique is discussed in Sect.~\ref{sec:error}. Subsequently, numerical examples and validation results are presented in Sect.~\ref{sec:NE}. Finally, we make the concluding remarks in Sect.~\ref{sec:con}.

\section{Numerical integration (PGQ)} \label{sec:NI}
\subsection{Computational procedure} \label{sec:proc}
There are two steps in the proposed numerical integration technique: (1) compute numerical integral in each integration cell; (2) assemble the matrices of integration cells to form the elemental matrices of the cut cell. FCM method also incorporate similar computational procedures. In FCM method~\cite{Duester2008Cmiamae}, the Jacobian terms are modified to establish relationship and map between the integration cell and the cut cell. On the other hand, in the proposed numerical integration technique, the FE formulation and Quadrature rule remain unchanged. This characteristics significantly improves its scalability and easy-to-implement to existing FEM solvers. The assembly procedure is computed through quadratic form transformation~\cite{gregory1981quadratic}. Therefore, this technique is termed as projection-based Gaussian quadrature (PGQ).  

The detailed algorithm of PGQ is demonstrated based on a general FE formulation below. Assuming the domain is spatially discretized by structured quadrilateral elements in Fig.~\ref{fig:demo}, the corresponding variational form of a general partial differential equation (PDE) over a cut cell can be defined as

\begin{eqnarray}
\mathcal{A}(\bm{v},\bm{d}_h) & = & \mathcal{L}(\bm{v}) \label{eq:elasticity}\\
\mathcal{A}(\bm{v},\bm{d}_h) & = & \int \limits_{\Omega} [\bm{L} \cdot \bm{v}]' \cdot \bm{D} \cdot [\bm{L} \cdot \bm{d}_h] d\Omega \nonumber \\
\mathcal{L}(\bm{v}) & = & \int \limits_{\Omega} [\bm{v}' \cdot \bm{b}] d\Omega + \int \limits_{\Gamma_H} [\bm{v}' \cdot \tilde{\bm{h}}_h] d\Gamma \nonumber
\end{eqnarray}
where $\mathcal{A}(\bm{v},\bm{d}_h)$ and $\mathcal{L}(\bm{v})$ are respectively bilinear and linear functionals. In Eq.~\eqref{eq:elasticity}, $\bm{v}$, $\bm{d}_h$, $\bm{D}$, $\bm{b}$ and $\tilde{\bm{h}}_h$ are test function vector, nodal value vector, coefficient matrix, volume source vector and prescribed traction vector respectively. $\bm{L}$ is denoted as a differential operator, where the prime symbol is a transpose operator. The strain matrix $\bm{B}$ is defined as $\bm{B} = \bm{L}\cdot\bm{N}$, where $\bm{N}$ is trial function matrix. In Bubnov-Galerkin method, the test function is chosen as trial function, $\bm{v} = \bm{N}$. Hence the elemental stiffness matrix and force matrix of the cut cell becomes,
\begin{eqnarray}
\bm{K}_c &=& \int \limits_{\Omega} [\bm{B}'_c \cdot \bm{D} \cdot \bm{B}_c] d\Omega \\
\bm{F}_c &=& \int \limits_{\Omega} [\bm{N}'_c \cdot \bm{b}] d\Omega + \int \limits_{\Gamma_H} [\bm{N}'_c \cdot \tilde{\bm{h}}_h] d\Gamma
\end{eqnarray} 
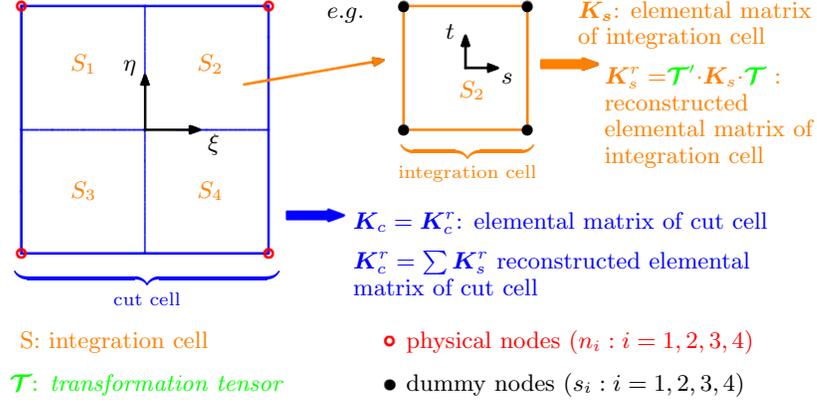
\begin{figure} \centering
	\hspace{-25pt}\includegraphics[trim=0.1cm 6.5cm 0.1cm 9.5cm,scale=0.4,clip]{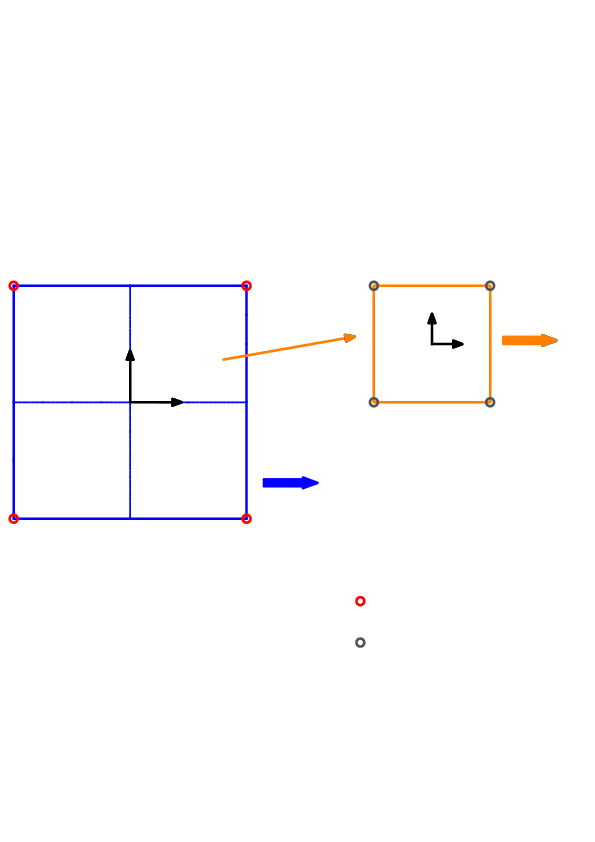}
	\begin{picture}(0,0)
	\put(-217,125){\small $\orange{S_1}$}
	\put(-169,125){\small $\orange{S_2}$}
	\put(-217,77){\small $\orange{S_3}$}
	\put(-169,77){\small $\orange{S_4}$}
	\put(-70,115){\small $\orange{S_2}$}
	\put(-165,95){\small $\xi$}
	\put(-197,125){\small $\eta$}
	\put(-54,120){\small $s$}
	\put(-75,137){\small $t$}
	\put(-120,145){\small $e.g.$}
	\put(-236,20){\small \orange{S: integration cell}}
	\put(-240,4){\small \green{$\bm{\mathcal{T}}$: \textit{transformation tensor}}}
	\put(-25,145){\small \orange{$\bm{K_s}$: elemental matrix}}
	\put(-25,135){\small \orange{of integration cell}}
	\put(-15,120){\small $\orange{\bm{K}^r_s =} \green{\bm{\mathcal{T}}'} \orange{\cdot\bm{K}_s \cdot} \green{\bm{\mathcal{T}}}$ \orange{:}}
	\put(-15,110){\small \orange{reconstructed}}
	\put(-15,100){\small \orange{elemental matrix of}}
	\put(-15,90){\small \orange{integration cell}}
	\put(-110,65){\small \blue{$\bm{K}_c = \bm{K}^r_c$: elemental matrix of cut cell}}
	\put(-110,50){\small \blue{$\bm{K}^r_c = \sum \bm{K}^r_s$ reconstructed elemental}}
	\put(-110,40){\small \blue{matrix of cut cell}}
	\put(-90,20){\small \red{physical nodes ($n_i: i=1,2,3,4$)}}
	\put(-90,4){\small dummy nodes ($s_i: i=1,2,3,4$)}
	\put(-96,6){\color{black}\circle*{4}}
	\put(-91,149.2){\color{black}\circle*{4}}
	\put(-91,102.5){\color{black}\circle*{4}}
	\put(-44.2,149.2){\color{black}\circle*{4}}
	\put(-44.2,102.5){\color{black}\circle*{4}}
	\put(-238,50){\blue{$\underbrace{\qquad \qquad \qquad \qquad \qquad}_{\text{cut cell}}$}}
	\put(-93,98){\orange{$\underbrace{\qquad \quad \qquad}_{\text{integration cell}}$}}
	\end{picture}
	\caption{Illustration of the general concept of PGQ in a bilinear quadrilateral element}
	\label{fig:demo}
\end{figure}
where the subscript $"c"$ refers to a cut cell. The standard Gauss quadrature rule is implemented in each integration cell with respect to its dummy nodes (black circle) in detailed view of Fig.~\ref{fig:demo}. The embedded interface in cut cell is assumed to align with the edges of integration cell. Therefore, in each integration cell, the variable values and their gradients on dummy nodes are approximated within a finite space of continuous function. 

Similar to the cut cell, the stiffness matrix and force vector of an integration cell are defined as
\begin{eqnarray}
\bm{K}_s &=& \int \limits_{\Omega} [\bm{B}'_s \cdot \bm{D} \cdot \bm{B}_s] d\Omega \label{eq:Ks} \\
\bm{F}_s &=& \int \limits_{\Omega} [\bm{N}'_s \cdot \bm{b}] d\Omega + \int \limits_{\Gamma_H} [\bm{N}'_s \cdot \tilde{\bm{h}}_h] d\Gamma \label{eq:Fs}
\end{eqnarray} 
where the subscript $"s"$ refers to an integration cell. The second term in Eq.~\eqref{eq:Fs} is Neumann boundary condition along the edge of an integration cell, which is associated with embedded interface in the cut cell.
As demonstrated in Fig.~\ref{fig:demo}, they are mapped via transformation tensor $\bm{\mathcal{T}}$ and assembled to form the reconstructed elemental matrices of a cut cell, e.g., $\bm{K}^r_c$, where the superscript $"r"$ denotes a reconstructed matrix. The transformation procedure is based on change of basis operation in a quadratic form. The elemental matrices are mapped between basis vectors of integration cell and cut cell.

\begin{algorithm}
	\caption{transformation-assembly} \label{co:algo1}
	\begin{algorithmic}[1]
		\State sub-divide cut cell into integration cells
		\For{i=no. of integration cell}
		\State construct $\mathcal{\bm{T}}$ for ith integration cell, Eq.~\ref{eq:T}
		\State Gaussian quadrature for ith integration cell
		\State $\bm{K}^r_{s}=\mathcal{\bm{T}}'
		\cdot \bm{K}_{s} \cdot \mathcal{\bm{T}}$, Eq.~\ref{eq:Kr} and~\ref{eq:Fr}
		\State sum up as $\bm{K}^r_c += \bm{K}^r_{s}$, Eq.~\ref{eq:KrAssem} and~\ref{eq:FrAssem}
		\EndFor
	\end{algorithmic}
\end{algorithm}
\begin{algorithm}
	\caption{assembly-transformation} \label{co:algo2}
	\begin{algorithmic}[1]
		\State sub-divide cut cell into integration cells
		\State construct $\mathcal{\bm{T}}^{as}$ (rectangular tensor) for all integration cells, Eq.~\ref{eq:T}
		\For{i=no. of integration cell}
		\State Gaussian quadrature for ith integration cell
		\State assembly as $\bm{K}^{as}_s$
		\EndFor
		\State $\bm{K}^r_c=\mathcal{\bm{T}}^{as'} 
		\cdot \bm{K}^{as}_s \cdot \mathcal{\bm{T}}^{as}$, Eq.~\ref{eq:KrAssem2} and~\ref{eq:FrAssem2}
	\end{algorithmic}
\end{algorithm}

Two types of computational sequence, Algorithm~\ref{co:algo1} and~\ref{co:algo2}, in PGQ are applicable, where the superscript $"as"$ denotes an assembled matrix based on standard assembly procedure in FEM. Both Algorithms are equivalent. Algorithm~\ref{co:algo1} is more computational efficient and preferred, because the operation of low-order matrix is involved. Nonetheless, Algorithm~\ref{co:algo2} demonstrates an important mathematical characteristics of PGQ, which will be discussed in Sect.~\ref{sec:Char}. In the next section, the construction of $\mathcal{\bm{T}}$ is discussed.

\subsection{Transformation tensor} \label{sec:trans}
\begin{figure} \centering
	\hspace{-25pt}\includegraphics[trim=0.01cm 2cm 1cm 6.2cm,scale=0.5,clip]{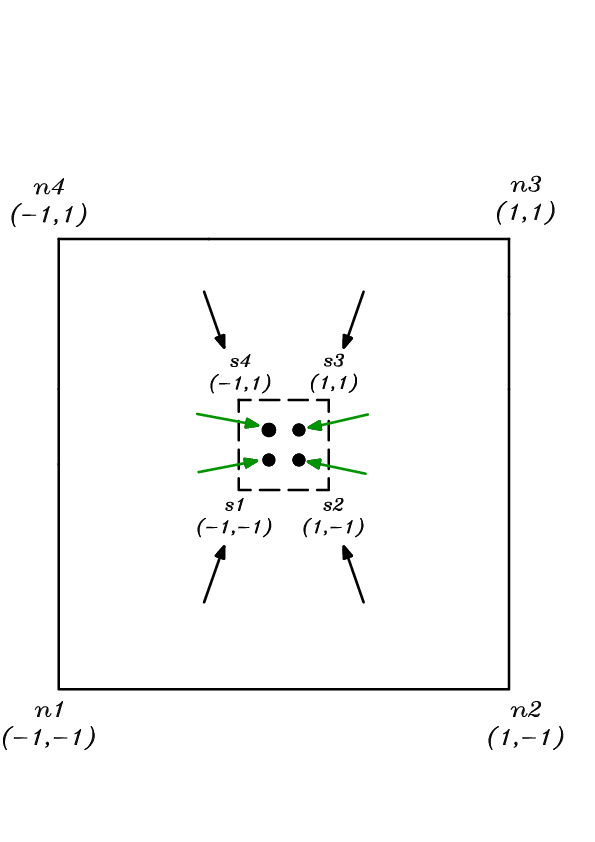}
	\begin{picture}(0,0)
	\put(-250,75){\small$d_s(s_1) = N_{c(i)}(s_1) d_c(n_i)$}
	\put(-138,75){\small$d_s(s_2) = N_{c(i)}(s_2) d_c(n_i)$}
	\put(-250,255){\small$d_s(s_4) = N_{c(i)}(s_4) d_c(n_i)$}
	\put(-138,255){\small$d_s(s_3) = N_{c(i)}(s_3) d_c(n_i)$}
	\put(-240,155){\small$V(\bm{x}_1) = $}
	\put(-252,143){\small$N_{s(i)}(\bm{x}_1) d_s(s_i)$}
	\put(-88,155){\small$V(\bm{x}_2) = $}
	\put(-100,143){\small$N_{s(i)}(\bm{x}_2) d_s(s_i)$}
	\put(-240,197){\small$V(\boldsymbol{x}_4) = $}
	\put(-252,184){\small$N_{s(i)}(\bm{x}_4) d_s(s_i)$}
	\put(-88,197){\small$V(\bm{x}_3) = $}
	\put(-100,184){\small$N_{s(i)}(\bm{x}_3) d_s(s_i)$}
	\end{picture}
	\caption{Transformation matrix: mapping procedure between integration cell and cut cell}
	\label{fig:ica}
\end{figure}
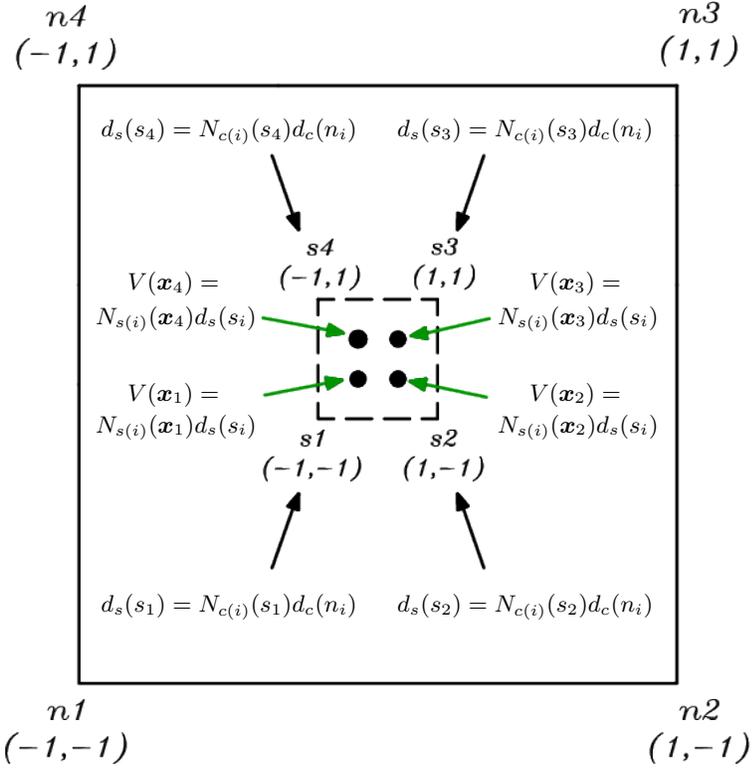
A detailed construction procedure of $\bm{\mathcal{T}}$ is demonstrated in Fig.~\ref{fig:ica}. A scalar value $V$ at Gauss points inside an integration cell, the solid-circles in Fig.~\ref{fig:ica}, are approximated from the dummy nodes, $s_i$. Simultaneously, the values on dummy nodes are approximated from the physical nodes, $n_i$. As a result, the scalar value on Gauss point $\bm{x}_k$, $V(\bm{x}_k)$, can be approximated from the physical nodes as shown in Eq.~\eqref{eq:V}..
\begin{eqnarray}
V(\bm{x}_k) & = & N_{s(j)}(\bm{x}_k)N_{c(i)}(s_j)d_c(n_i) \nonumber \\
& = & \tilde{N}_i(\bm{x}_k) d_c(n_i) \label{eq:V}
\end{eqnarray}  
where $\tilde{\bm{N}}(\bm{x})$ is a composed trial function vector. To put the aforementioned mapping procedure into a tensor form, a $\bm{\mathcal{T}}$ can be defined in Eq.~\eqref{eq:T}. The column $j$ of $\mathcal{\bm{T}}$ refers to the weights from a physical node $n_j$ to the dummy nodes $s_k$ of a cut cell. Therefore, $\tilde{\bm{N}}$ can be re-casted as the form in Eq.~\eqref{eq:Nt}. 
\begin{eqnarray}
\mathcal{T}_{kj} & = & N_{c(j)}(s_k) \label{eq:T} \\
\tilde{N}_{ij} & = & N_{s(k)}(\bm{x}_i) \mathcal{T}_{kj} \label{eq:Nt}
\end{eqnarray}
This transformation tensor is subsequently used to map the elemental matrices between the bases of integration cell and cut cell, as shown in Eq.~\eqref{eq:Kr} and~\eqref{eq:Fr}.
\begin{eqnarray}
K^r_{s(ij)} & = & \mathcal{T}_{ki} K_{s(kl)} \mathcal{T}_{lj} \label{eq:Kr}\\
F^r_{s(i)} & = & \mathcal{T}_{ki} F_{s(k)} \label{eq:Fr}
\end{eqnarray}
where $K^r_{s(ij)}$ and $F^r_{s(i)}$ are the reconstructed elemental matrices of an integration cell in component form.

Subsequently, the reconstructed elemental matrices of a cut cell is simply formed by a summation operation, as shown in Eq.~\eqref{eq:KrAssem} and~\ref{eq:FrAssem}. The $en$ parameter is the total number of integration cells in a cut cell.
\begin{eqnarray}
K^r_{c(ij)} & = & \sum \limits^{en}_{n=1} K^r_{s(ij)} (n) = \sum \limits^{en}_{n=1} \mathcal{T}_{ki}(n) K_{s(kl)}(n) \mathcal{T}_{lj}(n) \label{eq:KrAssem}\\
F^r_{c(i)} & = & \sum \limits^{en}_{n=1} F^r_{s(i)}(n) = \sum \limits^{en}_{n=1} \mathcal{T}_{ki}(n) F_{s(k)}(n) \label{eq:FrAssem}
\end{eqnarray}
The above demonstrates the computational procedure in Algorithm~\ref{co:algo1}.
As mentioned in Sect.~\ref{sec:proc}, the assembly procedure can be performed before transformation operation in Algorithm~\ref{co:algo2}. This assembly procedure the standard matrix assembly procedure in FEM. The transformation operation in Algorithm~\ref{co:algo2} is shown in Eq.~\eqref{eq:KrAssem2} and~\ref{eq:FrAssem2}. 
\begin{eqnarray}
K^r_{c(ij)} &=& \mathcal{T}^{as}_{ki} K^{as}_{s(kl)} \mathcal{T}^{as}_{lj} \label{eq:KrAssem2} \\
F^r_{c(i)} &=& \mathcal{T}^{as}_{ki} F^{as}_{s(k)} \label{eq:FrAssem2}
\end{eqnarray}
where $\mathcal{T}^{as}$ is a rectangular transformation tensor which has number of rows as $\bm{K}^{as}_s$ and number of columns as $\bm{K}^r_c$. The definition of $\mathcal{T}^{as}$ is identical with $\mathcal{T}$ in Eq.~\ref{eq:T}. In the next section, the implementation of PGQ is briefly discussed.
 
\subsection{Implementation of PGQ}
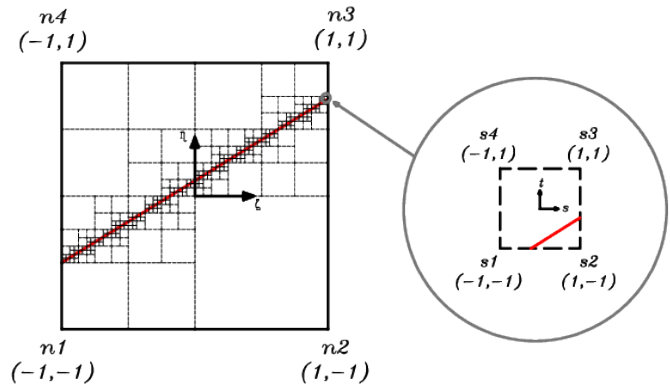
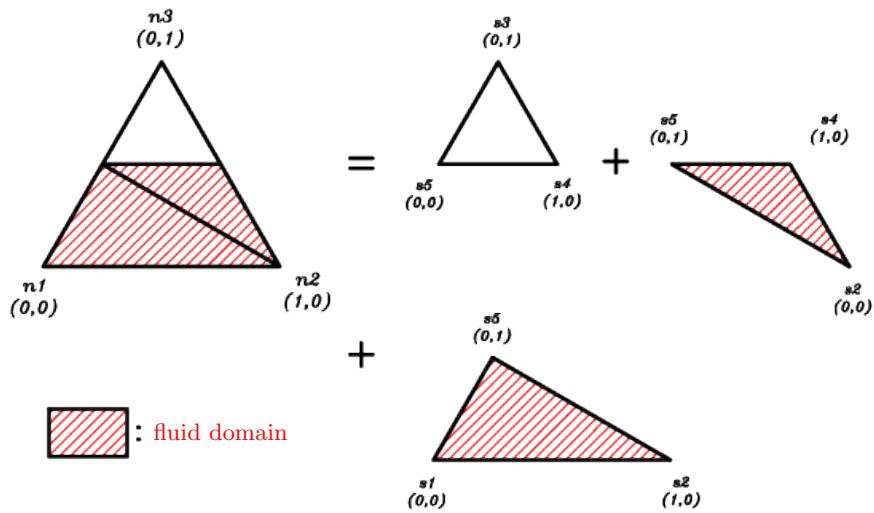
\begin{figure} \centering
	\begin{subfigure}[b]{1.0\textwidth}	
		\centering
		\hspace{-25pt}\includegraphics[trim=0.1cm 6cm 1.2cm 10.5cm,scale=0.45,clip]{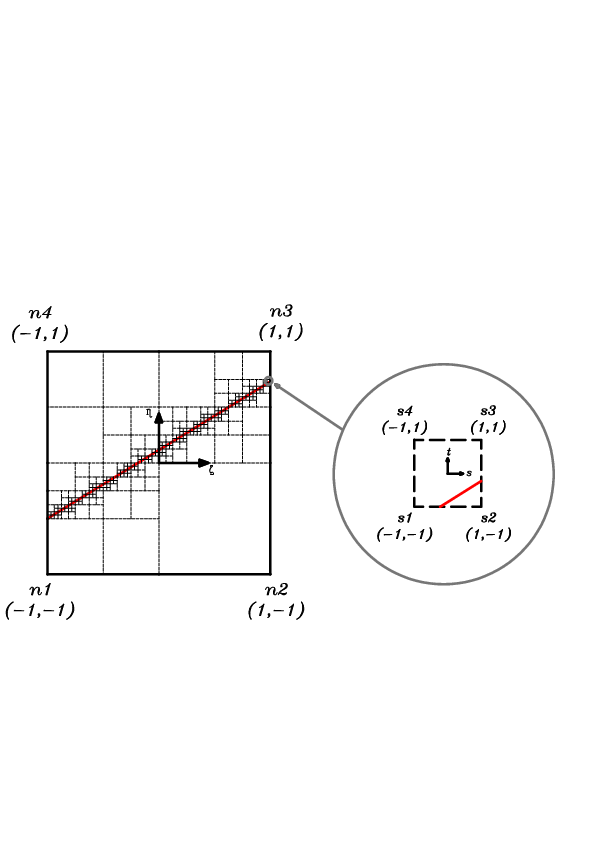}
		\caption{$\qquad$}
		\label{fig:quad}
	\end{subfigure}
	\begin{subfigure}[b]{1.0\textwidth}
		\centering
		\hspace{-25pt}\includegraphics[trim=0.1cm 7cm 0.1cm 9.5cm,scale=0.6,clip]{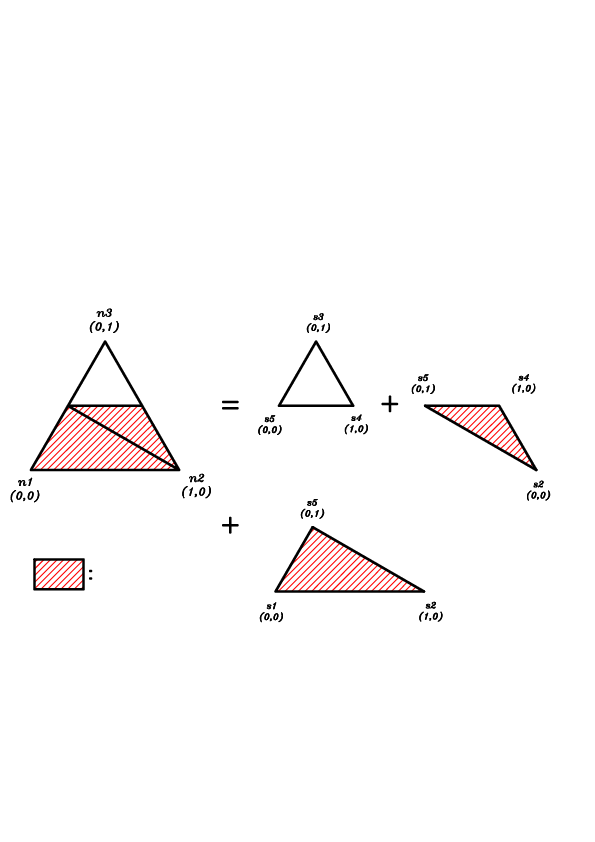}
		\begin{picture}(0,0)
		\put(-300,39){\small \red{fluid domain}}
		\end{picture}
		\caption{$\qquad$}
		\label{fig:tri}
	\end{subfigure}
	\caption{Implementation of PGQ: (a) quadtree refinement; (b) tessellation}
	\label{fig:quadTri}
\end{figure}

The proposed PGQ can be implemented via adaptive refinement in Fig.~\ref{fig:quad}, or tessellation in Fig.~\ref{fig:tri}. In quadtree adaptive refinement method, the mesh is locally refined along the embedded interface. It is able to capture embedded interface with strong geometric nonlinearity. However, its computational cost is relative high and the integration cells cannot accurately align with the embedded interface. On the other hand, tessellation method is more computational efficient. The tessellation method is well-established and considered as one of the standard numerical integration techniques in embedded interface problem, in which the integration cell is discretized such that its edges exactly align with the embedded interface. 

In quadtree adaptive refinement method, it is recommended to take Algorithm~\ref{co:algo1}, since it results into a huge number of integration cell. On the other hand, both Algorithm~\ref{co:algo1} and~\ref{co:algo2} can be efficiently implemented in tessellation method. In the next section, important characteristics of proposed PGQ will be discussed in detail.

\subsection{Characteristics of PGQ} \label{sec:Char}
In this section, important characteristics of proposed PGQ are discussed. They are (1) \emph{partition of unity} property, (2) recovery of Gauss quadrature (GQ), (3) projection in quadratic form and (4) reduced-order modeling.

The \emph{partition of unity} is one of the fundamental properties in FEM approximation. It can be simply proven that the composed trial function $\tilde{\bm{N}}$ in Sect.~\ref{sec:trans} do satisfy the \emph{partition of unity} property, as shown in Eq.~\eqref{eq:unity}. 
\begin{eqnarray}
\sum^4_{i=1} \tilde{N}_i(\alpha) & = & \sum \limits^4_{i=1} N_{j}(\alpha)N_{i}(\beta_j) = \sum^4_{i=1} N_{i}(\alpha) = 1 \label{eq:unity}
\end{eqnarray}

The second characteristics is the exact recovery of Gauss quadrature by reconstructed elemental matrices. The following mathematical derivation in Eq.~\eqref{eq:proof} shows that the reconstructed elemental matrices, e.g., $\bm{K}^r_c$, exactly recovers the the elemental matrices computed from the standard Gauss quadrature numerical integration of the cut cell, e.g., $\bm{K}_c$.  
\begin{eqnarray}
\bm{K}^r_c & = & \sum \limits^{en}_{n=1} \bm{\mathcal{T}}'(n) \cdot \bm{K}_s (n) \cdot \bm{\mathcal{T}}(n) \nonumber \\
& = & \sum \limits^{en}_{n=1} \bm{\mathcal{T}}'(n) \cdot \int \limits_{\Omega} [ \bm{B}'_s(n) \cdot \bm{D} \cdot \bm{B}_s(n)] d\Omega \cdot \bm{\mathcal{T}}(n) \nonumber \\
& = & \sum \limits^{en}_{n=1} \int \limits_{\Omega} [\tilde{\bm{N}}'(n) \cdot \bm{L}' \cdot \bm{D} \cdot \bm{L} \cdot \tilde{\bm{N}}(n)] d\Omega \nonumber \\
& = & \sum \limits^{en}_{n=1} \sum \limits^{gp}_{g=1} [\tilde{\bm{B}}'(n,g) \cdot \bm{D} \cdot \tilde{\bm{B}}(n,g) |\bm{J}(n,g)| W(g)] \nonumber \\
& = & \sum \limits^{en \cdot gp}_{k=1} [\tilde{\bm{B}}'(k) \cdot \bm{D} \cdot \tilde{\bm{B}}(k) |\bm{J}(k)| W(k)] \nonumber\\
& = & \int \limits_{\Omega} \tilde{\bm{B}}' \cdot \bm{D} \cdot \tilde{\bm{B}} d\Omega \nonumber \\
& = & \bm{K}_c \label{eq:proof}
\end{eqnarray}
where $\tilde{\bm{B}}(k) = \bm{L} \cdot \tilde{\bm{N}}(k)$, $\bm{J}(k)$ and $W(k)$ are composed strain matrix, Jacobian matrix and the Gauss integration weights. The value of $en \cdot gp$ are the total number of the Gauss integration points within a cut cell, in which $en$ and $gp$ are respectively the total number of integration cell and the number of Gauss points within each integration cell. Furthermore, the exact recovery of Gauss quadrature is qualitatively shown by a lid-driven cavity flow problem in Sect.~\ref{sec:error}. It highlights the robustness of the proposed PGQ. as it works together with Gauss quadrature. 

In addition, it is noticed that the composed strain matrix $\tilde{\bm{B}}$ is derived based on iso-parametric formulation over continuous function space in integration cell; whereas its basis vector set is chosen as those of its cut cell. It guarantees an accurate approximation of the gradients of variable within the integration cell from the physical nodes of cut cell. Therefore, the stresses along the embedded interface can be approximated as those within a standard element of FEM formulation. Similar to $\tilde{\bm{B}}$, the Jacobian matrix $\bm{J}$ is computed based on the iso-parametric mapping of integration cell too. If there was an approximation error induced by embedded interface in an infinitesimal integration cell, the influence of this error can be minimized owing to its small Jacobian value $|\bm{J}|$. This is particular true in quadtree adaptive refinement method by discretizing sufficient small integration cells along the embedded interface.  

The third characteristic of PGQ is about its quadratic form. It is known that a matrix, e.g., $\bm{K}_s$ in Eq.~\eqref{eq:basis}, can be mapped back to its own basis function space using its unit basis vectors, e.g., $\bm{e}_1 = [1, 0, 0]'$ in Cartesian coordinate system of $\mathbb{R}^3$. 
\begin{eqnarray}
K_{s(i,j)} = \bm{e}'_i \cdot \bm{K}_s \cdot \bm{e}_j \label{eq:basis}
\end{eqnarray}
Similarly, it can be projected to other basis function spaces, provided an appropriate transformation tensor is defined. In PGQ, $\bm{\mathcal{T}}$ is constructed based on its trial functions, linearly independent vector set in Eq.~\eqref{eq:T}, such that the nodal values and their residuals are re-projected in the basis function space of the cut cell.

In variational principle, to find a set of discrete solutions in FEM formulation, which minimizes the residual of equation system in an integral sense over a computational domain subjected to boundary conditions, can be treated as a quadratic optimization problem. Because the FEM formulation results into a symmetric matrix system~\footnote{In Navier-Stokes equation, the resultant stiffness matrix $\bm{K}$ can be subdivided into symmetric matrix blocks} and a symmetric matrix can always be transformed into a quadratic form, the proposed PGQ is mathematically-robust and well-suited for FEM formulation. It is consistent with the origin derivation of FEM theorem.

As shown in Algorithm~\ref{co:algo2}, the proposed PGQ can be deemed as a ROM technique. Recollecting Fig.~\ref{fig:tri}, a constant-strain triangular (CST) cut cell is discretized into three CST integration cells, and two additional DoFs, dummy nodes, are introduced. Therefore, the numerical integration of a cut cell with embedded discontinuity, 3 DoFs, is projected to a higher-dimension space, 5 DoFs, where the nonlinear problem maybe linearly separable. After the numerical integrations are performed in a higher-dimension space, the resultant matrix system is projected back to a lower-dimension space via an appropriate transformation tensor $\bm{\mathcal{T}}$ in quadratic form. When the matrix system is projected back to a lower-dimension space, the accuracy of results is subjected to a sufficient number of DoFs. For the case in Fig.~\ref{fig:tri}, because there are only 3 DoFs in CST cut cell, the embedded interface cannot be approximated accurately and result into a local smoothing of the embedded discontinuity in a cut cell.

In this section, the introduction of our proposed PGQ technique is completed. In the next section, we are going to present the governing equations which we are solving in our unfitted FSI solver together with implemented time integration and FSI coupling schemes. 

\section{Governing equations and boundary conditions} \label{sec:govern}
\subsection{Incompressible Navier-Stokes equations}
In our developed FSI solver, we are solving for a moving rigid body submerged in an incompressible Newtonian fluid. Therefore, the incompressible Navier-Stokes equations, Eq.~\eqref{eq:ns1} to~\eqref{eq:ns4}, are implemented.
\begin{eqnarray}
\rho^f \Big( \frac{\partial \bm{u}^f}{\partial t} + \bm{u}^f \cdot \nabla \bm{u}^f \Big)- \nabla \cdot \boldsymbol{\sigma}\{\bm{u}^f,p\} & = & \rho^f \bm{g}^f \quad \forall \bm{x} \in \bm{\Omega}^f (t) \label{eq:ns1}\\
\nabla \cdot \bm{u}^f & = & 0 \qquad\;\; \forall \bm{x} \in \bm{\Omega}^f(t) \label{eq:ns2}\\
\bm{u}^f & =& \tilde{\bm{u}}^f \quad\;\;\; \forall \bm{x} \in \bm{\Gamma}^f_D(t)\label{eq:ns3}\\
\boldsymbol{\sigma}\{\bm{u}^f,p\}\cdot \bm{n}^f & = & \tilde{\bm{h}}^f \quad\;\;\; \forall \bm{x} \in \bm{\Gamma}^f_H(t) \\
\bm{u}^f &=& \bm{u}^f_0 \quad\;\;\; \forall \bm{x} \in \bm{\Omega}^f(0) \label{eq:ns4}
\end{eqnarray}
where $\rho^f$, $\bm{u}^f$, $\bm{u}^f_0$, $\bm{g}^f$, $\tilde{\bm{u}}^f$, $\tilde{\bm{h}}^f$ and $\bm{n}^f$ are respectively the fluid density, fluid velocity vector, initial fluid velocity vector, fluid unit body force vector, prescribed fluid velocity, prescribed fluid traction and outward normal vector of fluid domain. The supperscript $f$ and $s$ refer to fluid and solid respectively. The spatial domain, Dirichlet and Neumann boundaries are respectively denoted as $\bm{\Omega}$, $\bm{\Gamma}_D$ and $\bm{\Gamma}_H$, where $\bm{\Gamma}_D$ and $\bm{\Gamma}_H$ are complementary subsets of $\bm{\Gamma}$, $\bm{\Gamma} = \bm{\Gamma}_D \cup \bm{\Gamma}_H$ and $\bm{\Gamma}_D \cap \bm{\Gamma}_H = \O$. The Dirichlet and Neumann boundary conditions respectively are imposed along $\bm{\Gamma}_D$ and $\bm{\Gamma}_H$ as shown below.
\begin{eqnarray}
\bm{u}^f & = & \tilde{\bm{u}}^f \quad \forall \bm{x} \in \bm{\Gamma}^f_D(t) \label{eq:EBC}\\
\bm{h}^f  & = & \tilde{\bm{h}}^f \quad\; \forall \bm{x} \in \bm{\Gamma}^f_H(t) \label{eq:NBC}
\end{eqnarray}
where $\bm{h} = \bm{\sigma}\cdot\bm{n}$ refers to the surface stresses. $\bm{\sigma}$ is the Cauchy stress tensor and defined as 
\begin{eqnarray}
\bm{\sigma}\{\bm{u}^f,p\} & = & -p \bm{I} + 2 \mu D(\bm{u}^f) \\ 
D(\bm{u}^f) & = & \frac{1}{2} \Big[ \nabla \bm{u}^f + (\nabla \bm{u}^f)' \Big]
\end{eqnarray}
The stress tensor is written as the summation of its isotropic and deviatoric tensor ($D(\bm{u}^f)$) parts. Here, $p$, $\mu$ and $\bm{I}$ refer to the fluid pressure, dynamic viscosity and identity matrix respectively. 

\subsection{Rigid-body dynamics}
The equations governing dynamics of a rigid body is simply implemented as Eq.~\eqref{eq:struct}.
\begin{eqnarray}
m^s \bm{a}^s + \bm{c}^s \cdot \bm{u}^s + \bm{k}^s \cdot \bm{d}^s & = & \bm{h}^s \quad \forall \bm{x} \in \bm{\Omega}^s(t) \label{eq:struct} \\
\nonumber \bm{a}^s = \frac{\partial^2 \bm{d}^s}{\partial t^2}; \quad \bm{u}^s &=& \frac{\partial \bm{d}^s}{\partial t}\\
\nonumber \bm{c}^s = 2\xi \sqrt{\bm{k}^s m^s}; \quad \bm{k}^s & = & 4 \pi^2 \bm{f}^2_n m^s \\
\nonumber U_r = U/(f_{ny} D); \quad m^s & = & m^* (0.25\pi D^2 L \rho^f)
\end{eqnarray}
where $\xi$, $m^*$, $\bm{f}_n = [f_{nx}, f_{ny}]'$, $D$ and $L$ are the damping ratio, mass ratio, structural frequency vector, diameter of cylinder and spanwise length of cylinder respectively. $\bm{c}^s$ and $\bm{k}^s$ refer to damping and stiffness coefficients respectively. The reduced velocity of the cylinder, $U_r$, is based on the structural frequency in the transverse direction, $f_{ny}$. In the present formulation, it is assumed that the structural frequencies in transverse and streamwise direction are identical, $f_{nx}/f_{ny} = 1.0$. $\bm{h}^s = [h^s_x, h^s_y]'$ represents the external force exerted on the cylinder surface. In FSI problems, these external forces are hydrodynamic forces exerted by fluid around the surface of cylinder.

\subsection{Interface constraints and Fluid-structure interaction}
To couple fluid and structure, velocity and traction constraints should be satisfied. The velocity constraint requires the fluid and structure interfaces align and move at the same velocity, as shown in Eq.~\eqref{eq:v_con}. On the other hand, the equilibirum of stresses, Eq.~\eqref{eq:t_con}, has to be enforced along the fluid-structure interface to satisfy the traction constraint, where $\bm{n}^f = -\bm{n}^s$ and the superscript $"fs"$ refers to fluid-structure interface.
\begin{eqnarray}
\bm{u}^f &=& \bm{u}^s \quad \forall \bm{x} \in \bm{\Gamma}^{fs}(t) \label{eq:v_con}\\
\bm{\sigma}^f \cdot \bm{n}^f + \bm{\sigma}^s \cdot \bm{n}^s &=& 0 \quad\;\; \forall \bm{x} \in \bm{\Gamma}^{fs}(t) \label{eq:t_con} \\
\implies \bm{h}^f + \bm{h}^s &=& 0 \quad\;\; \forall \bm{x} \in \bm{\Gamma}^{fs}(t) \nonumber
\end{eqnarray}

The fluid and structural governing equations can be coupled in either monolithic or staggered-partitioned scheme. In monolithic/fully-implicit scheme~\cite{Blom1998Cmiamae}, the overall system equation consists of the variables of fluid and structure. They are solved indiscriminantly and simultaneously. The monolithic formulation is robust, stable at relative large time steps, and its solution converges rapidly. Albeit monolithic schemes have the energy conservation property, their computational cost is high and typically require a significant recast in both existing fluid and structural solvers. On the other hand, the staggered-partioned scheme can be conveniently implemented to existing fluid and structural solvers. The staggered-partitioned schemes can be further classified into strongly-coupled~\cite{Dettmer2006CMiAMaE,Jaiman2016CFb,Kadapa2016CMiAMaE} or weakly-coupled schemes~\cite{Dettmer2013IJfNMiE,Placzek2009CF}. In this work, a staggered-partitioned, weakly-coupled and second-order accurate scheme~\cite{Dettmer2013IJfNMiE} is implemented. Please refer to~\cite{Dettmer2013IJfNMiE} for detailed algorithm.

\subsection{Integration in time}
To deal with moving embedded boundaries in FSI problems, the second-order accurate and unconditional stable generalized-$\alpha$ method~\cite{Chung1993Joam} and~\cite{Jansen2000Cmiamae} are implemented in time integration for structural equation and Navier-Stokes equations respectively. The detailed formulation for structural equation can be summarized as,

\begin{eqnarray}
\mathbf{d}^s_{n+1} & = & \mathbf{d}^s_{n} + \Delta t \mathbf{u}^s_n + \Delta t^2\big(\big(\frac{1}{2} - \beta^s\big)\mathbf{a}^s_{n} + \beta^s \mathbf{a}^s_{n+1}\big) \label{eq:ga1}\\
\mathbf{u}^s_{n+1} & = & \mathbf{u}^s_n + \Delta t \big((1-\gamma^s)\mathbf{a}^s_n + \gamma^s \mathbf{a}^s_{n+1}\big) \label{eq:ga2}\\
\mathbf{d}^s_{n+\alpha^s_f} &= & (1-\alpha^s_f) \mathbf{d}^s_n + \alpha^s_f \mathbf{d}^s_{n+1} \\
\mathbf{u}^s_{n+\alpha^s_f} & = & (1-\alpha^s_f) \mathbf{u}^s_n + \alpha^s_f \mathbf{u}^s_{n+1} \\
\mathbf{a}^s_{n+\alpha^s_m} & = & (1-\alpha^s_m) \mathbf{a}^s_n + \alpha^s_m \mathbf{a}^s_{n+1} \\
\mathbf{F}^s_{n+\alpha^s_f} & = & (1-\alpha^s_f) \mathbf{F}^s_n + \alpha^s_f \mathbf{F}^s_{n+1}
\end{eqnarray}
where $\mathbf{d}^s_n$, $\mathbf{u}^s_n$ and $\mathbf{a}^s_n$ refer to the displacement, velocity and acceleration of cylinder at time $t=n$.
where $\alpha^s_m$, $\alpha^s_f$, $\gamma^s$ and $\beta^s$ are defined by~\cite{Chung1993Joam} as 

\begin{eqnarray}
\alpha^s_m & := & \frac{2-\rho^s_{\infty}}{\rho^s_{\infty}+1}; \qquad \alpha^s_f := \frac{1}{\rho^s_{\infty} + 1} \nonumber\\
\gamma^s  &:=& 0.5 + \alpha^s_m - \alpha^s_f; \qquad \beta^s := 0.25 (1+\alpha^s_m - \alpha^s_f)^2 
\end{eqnarray}



Similarly, the generalized-$\alpha$ method for Navier-Stokes equations is listed below,

\begin{eqnarray}
\bm{u}^f_{n+1} &=& \bm{u}^f_n + \Delta t \big[(1-\gamma^f)\frac{\partial \bm{u}^f_n}{\partial t} + \gamma^f \frac{\partial \bm{u}^f_{n+1}}{\partial t}\big] \\
\bm{u}^f_{n+\alpha^f_f} &=&(1-\alpha^f_f) \bm{u}^f_n + \alpha^f_f \bm{u}^f_{n+1}\\
\frac{\partial \bm{u}^f_{n+\alpha^f_m}}{\partial t} &=& (1-\alpha^f_m) \frac{\partial \bm{u}^f_n}{\partial t} + \alpha^f_m \frac{\partial \bm{u}^f_{n+1}}{\partial t}\\
\alpha^f_m &:=& 0.5\frac{3-\rho^f_{\infty}}{1+\rho^f_{\infty}}; \quad \alpha^f_f := \frac{1}{1+\rho^f_{\infty}}; \quad \gamma^f := 0.5+\alpha^f_m -\alpha^f_f
\end{eqnarray}
Here $\rho^s_{\infty} \in$ [0,1] and $\rho^f_{\infty} \in$ [0,1] respectively are the spectral radius, which control the amount of numerical high-frequency damping in the temporal schemes. In this work. $\rho^s_{\infty} = \rho^f_{\infty} = 0.2$ are chosen for all numerical simulations.


\section{Variational form of unfitted stabilized finite element formulation}\label{sec:formulation}
The complete stabilized FE formulation of Navier-Stokes equations with an embedded interface is summarized in Eq.~\ref{eq:all}, where $\mathcal{A}^{G}([\bm{v}^f,q],[\bm{u}^f_h,p_h])$ and $\mathcal{L}^{G}([\bm{v}^f,q],[\bm{u}^f_h,p_h])$ are the bilinear and linear forms derived from classical Galerkin method. $\mathcal{A}^{S}([\bm{v}^f,q],[\bm{u}^f_h,p_h])$ attributes to Pretrov-Galerkin formulation, which enables equal approximation function spaces between velocity and pressure. $\mathcal{A}^{N}([\bm{v}^f,q],[\bm{u}^f_h,p_h])$ is the terms of Nitsche's method for weakly imposing Dirichlet boundary condition along an embedded interface. In addition, $\mathcal{A}^{GP}([\bm{v}^f,q],[\bm{u}^f_h,p_h])$ is the ghost penalty terms to optimize the jump of quantities across edges in the cut cell.
\begin{eqnarray}
\mathcal{A}^{G}([\bm{v}^f,q],[\bm{u}^f_h,p_h]) + \mathcal{A}^{S}([\bm{v}^f,q],[\bm{u}^f_h,p_h]) + \mathcal{A}^{N}([\bm{v}^f,q],[\bm{u}^f_h,p_h]) & & \nonumber\\
+ \mathcal{A}^{GP}([\bm{v}^f,q],[\bm{u}^f_h,p_h]) = \mathcal{L}^{G}([\bm{v}^f,q],[\bm{u}^f_h,p_h])& & \label{eq:all}
\end{eqnarray}
The detailed formulations of the terms in Eq.~\ref{eq:all} are presented in the following sections.

\subsection{Stabilized variational form of Navier-Stokes equations}
The variational form of incompressible Navier-Stokes, Eq.~\eqref{eq:ns1} and~\eqref{eq:ns2}, based on classical Galerkin formulation is in Eq.~\eqref{eq:Gal}.
\begin{eqnarray}
\mathcal{A}^{G}([\bm{v}^f,q],[\bm{u}^f_h,p_h]) = \mathcal{L}^{G}([\bm{v}^f,q],[\bm{u}^f_h,p_h]) & & \nonumber\\
\implies \int \limits_{\Omega^f(t)}\bm{v}^f \Big[\rho^f \Big(\frac{\partial \bm{u}^f_h}{\partial t} + (\bm{u}^f_h \cdot \nabla) \bm{u}^f_h \Big) - \nabla \cdot \bm{\sigma} \{\bm{u}^f_h,p_h\} -\rho^f \bm{g}^f \Big]d\Omega & &\nonumber \\
+ \int \limits_{\Omega^f(t)} q [\nabla \cdot \bm{u}^f_h] d \Omega & & \nonumber\\ 
= \int \limits_{\Gamma^f_H(t)}\bm{v}^f \cdot \tilde{\bm{h}}^f_h d\Gamma \qquad \forall [\bm{v}^f,q] \in \hat{\bm{\mathcal{V}}}_h \times \hat{\mathcal{Q}}_h \subset \hat{\bm{{\mathcal{V}}}} \times \hat{\mathcal{Q}} & & \label{eq:Gal} 
\end{eqnarray}
where $[\bm{v}^f,q]'$ is the vector of test functions for velocity and pressure of fluid. The vector-valued trial and test function spaces $\boldsymbol{\mathcal{V}}$ and $\hat{\boldsymbol{\mathcal{V}}}$ of velocity are defined as 
\begin{eqnarray}
\boldsymbol{\mathcal{V}} &=& \{\bm{v}^f \in \bm{\mathcal{H}}^1 (\bm{\Omega}^f(t)): \bm{v}^f=\tilde{\bm{v}}^f \quad \forall \bm{x} \in \bm{\Gamma}^f_D(t)\} \nonumber\\
\hat{\boldsymbol{\mathcal{V}}} &=& \{\bm{v}^f  \in \bm{\mathcal{H}}^1 (\bm{\Omega}^f(t)): \bm{v}^f =\mathbf{0} \quad\;\; \forall \bm{x} \in \bm{\Gamma}^f_D(t) \}
\end{eqnarray} 
On the other hand, the scalar-valued trial and test function spaces $\mathcal{Q}$ and $\hat{\mathcal{Q}}$ of pressure are defined as 
\begin{eqnarray}
\mathcal{Q} &=& \{q \in \mathcal{H}^1 (\bm{\Omega}^f(t)): q=\tilde{p} \quad \forall \bm{x} \in \bm{\Gamma}^f_D(t)\} \nonumber\\
\hat{\mathcal{Q}} &=& \{q \in \mathcal{H}^1 (\bm{\Omega}^f(t)): q=0 \quad \forall \bm{x} \in \bm{\Gamma}^f_D(t) \}
\end{eqnarray} 
where $\boldsymbol{\mathcal{H}}^1$ refers to the Sobolev space, in which $[(\bm{v}^f)^2,q^2]$ and  $[|\nabla\bm{v}^f|^2,|\nabla q|^2]$ have finite integrals within $\bm{\Omega}^f(t)$ and allows discontinuous derivatives. Their corresponding discrete function spaces are denoted with subscript $"h"$, e.g., $\hat{\mathcal{Q}}_h$. 
In this work, a residual-based stabilization technique, Petrov-Galerkin method~\cite{Brooks1982Cmiamae,Shakib1991CMiAMaE,Tezduyar1992CMiAMaEa,Franca1992CMiAMaE}, in Eq.~\eqref{eq:stab} is implemented to ensure the residual of equation system is minimized in a (weak) integral sense over each element.
Here $\bm{G}$ and $C_I$ are respectively element cotravariant metric tensor and a positive constant independent upon mesh size~\cite{harari1992c}.

\begin{eqnarray}
\mathcal{A}^{S}([\bm{v}^f,q],[\bm{u}^f_h,p_h]) = \sum \limits^{n_{el}}_{e=1} \int \limits_{\Omega^f(t)} \tau_m \Big[\rho^f (\bm{u}^f_h \cdot \nabla) \bm{v}^f -\mu \nabla^2 \bm{v}^f + \nabla q \Big]
\cdot \Big[ \rho^f(\frac{\partial \bm{u}^f_h}{\partial t}  & & \nonumber\\
+ {(\bm{u}^f_h \cdot \nabla})\bm{u}^f_h-\bm{g}^f) - \mu\nabla^2 \bm{u}^f_h + \nabla p_h \Big] d\Omega + \sum \limits^{n_{el}}_{e=1} \int \limits_{\Omega^f(t)} \tau_c \rho^f (\nabla \cdot \bm{v}^f)(\nabla \cdot \bm{u}^f_h) d\Omega & & \label{eq:stab} \\ \qquad \forall [\bm{v}^f,q] \in \hat{\bm{\mathcal{V}}}_h \times \hat{\mathcal{Q}}_h \subset \hat{\bm{\mathcal{V}}} \times \hat{\mathcal{Q}} \nonumber \\
\nonumber \tau_m = \big[\big(\frac{2\rho^f}{\Delta t}\big)^2 + (\rho^f)^2 \bm{u}^f_h \cdot \bm{G} \bm{u}^f_h + C_I (\mu)^2 \bm{G}:\bm{G}\big]^{-0.5} & & \\
\nonumber \tau_c = (tr(\bm{G})\tau_m)^{-1}; \quad \bm{G}  =  \frac{\partial \bm{\xi}'}{\partial \bm{x}} \frac{\partial \bm{\xi}}{\partial \bm{x}} & &
\end{eqnarray}

\subsection{Nitsche's method}
To weakly impose Dirichlet boundary condition along the embedded interface, Nitsche's method is implemented. The terms of Nitsche's method are shown in Eq.~\eqref{eq:nitsche}.
\begin{eqnarray}
\mathcal{A}^{N}([\bm{v}^f,q],[\bm{u}^f_h,p_h]) = \gamma_{1} \int \limits_{\Gamma^{fs}(t)} \bm{v}^f \cdot (\bm{u}^f_h-\tilde{\bm{u}}^f_h) d\Gamma & &\nonumber \\
- \int \limits_{\Gamma^{fs}(t)} \bm{v}^f \cdot (\bm{\sigma}\{\bm{u}^f_h,p_h\} \cdot \bm{n}^f) d\Gamma 
- \gamma_{2} \int \limits_{\Gamma^{fs}(t)} (\bm{\sigma}\{\bm{v}^f,q\} \cdot \bm{n}^f) \cdot (\bm{u}^f_h-\tilde{\bm{u}}^f_h) d\Gamma & & \label{eq:nitsche} \\
\qquad \forall [\bm{v}^f,q] \in \hat{\bm{\mathcal{V}}}_h \times \hat{\mathcal{Q}}_h \subset \hat{\bm{\mathcal{V}}} \times \hat{\mathcal{Q}} & & \nonumber
\end{eqnarray}
Either symmetric-variant $\gamma_2 = 1$ or unsymmetric-variant $\gamma_2 = -1$ can be implemented. The penalty term is chosen within an appropriate range $\gamma_1 \in [\mu \frac{10^2}{L},$ $\mu\frac{10^3}{L}]$ \cite{Benk2012} for symmetric-variant, where $L$ is the characteristic element length, or $\gamma_1 = 0.0$ for unsymmetric-variant~\cite{Burman2012SJoNA}. As the solutions proceed to convergence $\bm{u}^f_h \approx \tilde{\bm{u}}^f_h$, the first and third penalty terms vanish. 

\subsection{Ghost Penalty Method}
The cut cell is separated by an embedded interface, e.g., the blue circle in Fig.~\ref{fig:schem}, into a fluid domain and a fictitious domain respectively. If the physical part is very small, some basis functions have little support inside the physical domain. It leads to large system matrix condition numbers.
The ghost penalty method~\cite{Burman2010CRM} is implemented along the edges of cut cells, the red edges in Fig.~\ref{fig:schem}, to alleviate the jumps of quantities in cut cells. 
A comprehensive study of the performance of ghost penalty terms was reported by Dettmer et al, 2016~\cite{Dettmer2016CMiAMaE}. The specific terms are listed in Eq.~\eqref{eq:gp}. 
\begin{figure} \centering
	\hspace{-25pt}\includegraphics[trim=0.1cm 7cm 0.1cm 10cm,scale=0.45,clip]{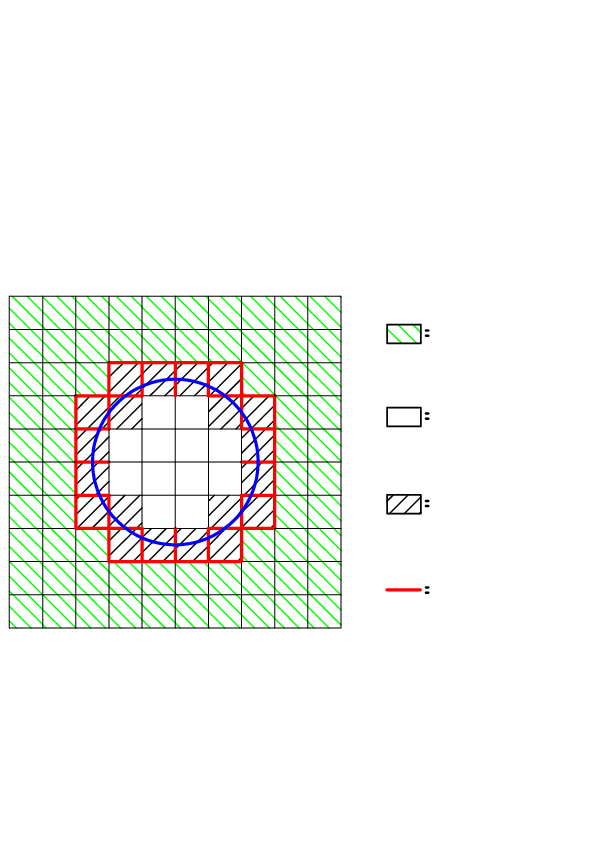}
	\begin{picture}(0,0)
	\put(-70,139){\small fluid domain}
	\put(-70,101){\small fictitious domain}
	\put(-70,61){\small cut cell}
	\put(-70,22){\small edges with ghost penalty}
	\end{picture}
	\caption{Schematic diagrams of embedded interface with ghost penalty terms}
	\label{fig:schem}
\end{figure}
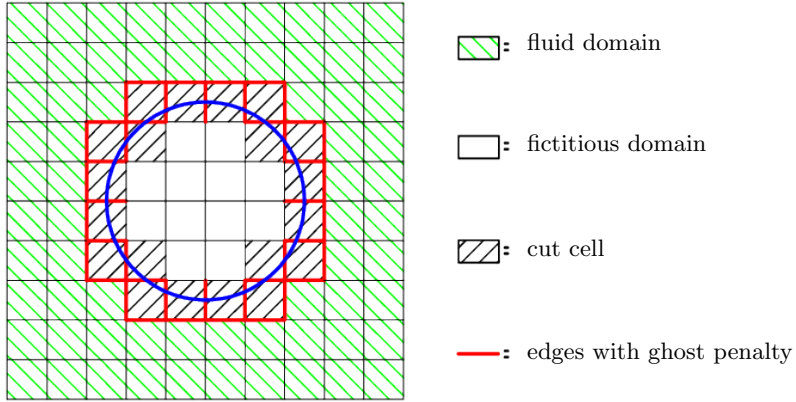
\begin{eqnarray}
\mathcal{A}^{GP}([\bm{v}^f,q],[\bm{u}^f_h,p_h]) = \beta^u_{gp} \mu G_1 (\bm{v}^f,\bm{u}^f_h) + \beta^p_{gp} \mu^{-1} g_3 (q,p_h) & &\label{eq:gp}\\
\quad \forall [\bm{v}^f,q] \in \hat{\bm{\mathcal{V}}}_h \times \hat{\mathcal{Q}}_h \subset \hat{\bm{\mathcal{V}}} \times \hat{\mathcal{Q}} & &\nonumber\\
g_{\phi}(q,p_h) = \sum \limits^e_{k=1} l^{2(\alpha-1)+\phi}_k \int \limits_{\Gamma^f(t)} [[\frac{\partial^{\alpha}q}{\partial \bm{n}^{f(\alpha)}}]][[\frac{\partial^{\alpha} p_h}{\partial \bm{n}^{f(\alpha)}}]] dl_k & &\nonumber \\
\nonumber G_{\phi}(\bm{v}^f,\bm{u}^f_h)  = \sum \limits^e_{k=1} \sum \limits^d_{i=1} l^{2(\alpha-1)+\phi}_k \int \limits_{\Gamma^f(t)} [[\frac{\partial^{\alpha} v^f_{(i)}}{\partial \bm{n}^{f(\alpha)}}]][[\frac{\partial^{\alpha}u^f_{h(i)}}{\partial \bm{n}^{f(\alpha)}}]]dl_k & &\\ \nonumber
\end{eqnarray}
where the penalty parameters is chosen as $\beta^u_{gp} = \beta^p_{gp} = 0.02$~\cite{Dettmer2016CMiAMaE} for simulations in this work. The superscripts $"u"$ and $"p"$ respectively refer to velocity and pressure. The subscript $"gp"$ shows these terms attribute to ghost penalty terms. $e$ and $d$ respectively are number edges of cut cell imposed with ghost penalty terms and dimension of problem. $\alpha$, $\phi$ and $l$ are order of derivative, the notation parameter and element characteristic length respectively. $[[\cdot]]$ denotes the jump of quantity across the element edge.

Therefore, the overall numerical formulation of Navier-Stokes equations with embedded interface is summarized as, 
\begin{eqnarray}
\int \limits_{\Omega^f(t)}\bm{v}^f \Big[\rho^f \Big(\frac{\partial \bm{u}^f_h}{\partial t} + (\bm{u}^f_h \cdot \nabla) \bm{u}^f_h \Big) - \nabla \cdot \bm{\sigma} \{\bm{u}^f_h,p_h\} - \rho^f \bm{g}^f \Big]d\Omega & &\nonumber\\
+ \int \limits_{\Omega^f(t)} q [\nabla \cdot \bm{u}^f_h] d \bm{\Omega}
+ \sum \limits^{n_{el}}_{e=1} \int \limits_{\Omega^f(t)}\tau_m \Big[ \rho^f (\bm{u}^f_h \cdot \nabla) \bm{v}^f -\mu \nabla^2 \bm{v}^f + \nabla q \Big] & & \nonumber\\
\cdot \Big[ \rho^f(\frac{\partial \bm{u}^f_h}{\partial t} + {(\bm{u}^f_h \cdot \nabla})\bm{u}^f_h-\bm{g}^f)
- \mu \nabla^2 \bm{u}^f_h + \nabla p_h \Big] d\Omega & &\nonumber\\
 + \sum \limits^{n_{el}}_{e=1} \int \limits_{\Omega^f(t)} \tau_c \rho^f (\nabla \cdot \bm{v}^f)(\nabla \cdot \bm{u}^f_h) d\Omega
+\gamma_{1} \int \limits_{\Gamma^{fs}(t)} \bm{v}^f \cdot (\bm{u}^f_h-\tilde{\bm{u}}^f_h) d\Gamma & & \nonumber\\
 - \int \limits_{\Gamma^{fs}(t)} \bm{v}^f \cdot (\bm{\sigma}\{\bm{u}^f_h,p_h\} \cdot \bm{n}^f) d\Gamma
-\gamma_{2} \int \limits_{\Gamma^{fs}(t)} (\bm{\sigma}\{\bm{v}^f,q\} \cdot \bm{n}^f) \cdot (\bm{u}^f_h-\tilde{\bm{u}}^f_h) d\Gamma & & \nonumber\\
+\sum \limits^e_{k=1} l^{2(\alpha-1)+\phi}_k \int \limits_{\Gamma^f(t)} [[\frac{\partial^{\alpha}q}{\partial \bm{n}^{f(\alpha)}}]][[\frac{\partial^{\alpha} p_h}{\partial \bm{n}^{f(\alpha)}}]] dl_k & &\nonumber\\
+ \sum \limits^e_{k=1} \sum \limits^d_{i=1} l^{2(\alpha-1)+\phi}_k \int \limits_{\Gamma^f(t)} [[\frac{\partial^{\alpha} v^f_{(i)}}{\partial \bm{n}^{f(\alpha)}}]][[\frac{\partial^{\alpha}u^f_{h(i)}}{\partial \bm{n}^{f(\alpha)}}]]dl_k & &\nonumber\\
= \int \limits_{\Gamma^f_H(t)}\bm{v}^f \cdot \tilde{\bm{h}}^f_h d\Gamma \qquad \forall [\bm{v}^f,q] \in \hat{\bm{\mathcal{V}}}_h \times \hat{\mathcal{Q}}_h \subset \hat{\bm{\mathcal{V}}} \times \hat{\mathcal{Q}} & &
\end{eqnarray}

\section{Convergence analysis} \label{sec:error}
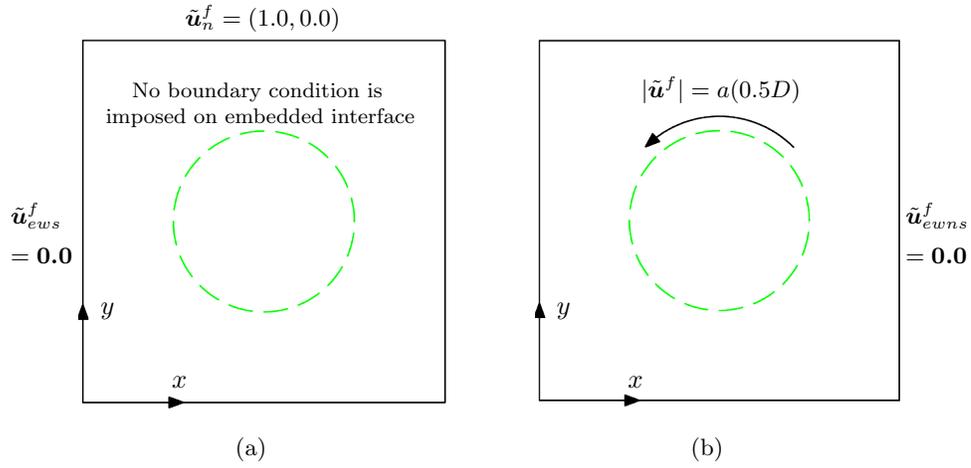
\begin{figure} \centering
	\begin{subfigure}[b]{0.5\textwidth}
		\centering
		\hspace{-25pt}\includegraphics[trim=0.3cm 3.5cm 1cm 6cm,scale=0.253,clip]{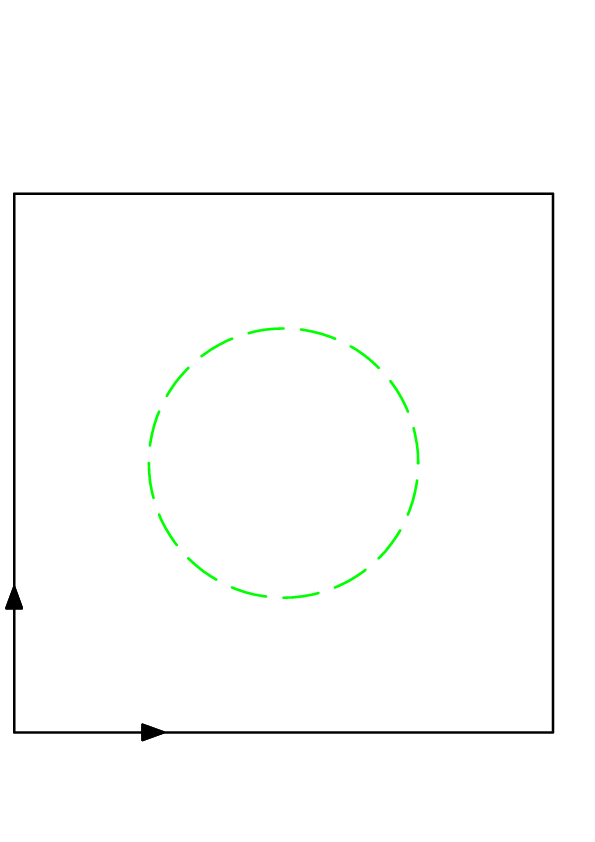}
		\begin{picture}(0,0)
		\put(-125,118){\footnotesize No boundary condition is}
		\put(-135,108){\footnotesize imposed on embedded interface}
		\put(-171,70){\small $\tilde{\bm{u}}^f_{ews}$}
		\put(-171,55){\small $ = \bm{0.0}$}
		\put(-105,145){\small $\tilde{\bm{u}}^f_n = (1.0, 0.0)$}
		\put(-137,36){$y$}
		\put(-110,8){$x$}
		\end{picture}
		\caption{$\qquad \qquad$}
		\label{fig:schLid}
	\end{subfigure}%
	\begin{subfigure}[b]{0.5\textwidth}	
		\centering
		\hspace{-25pt}\includegraphics[trim=0.3cm 3.5cm 1cm 6cm,scale=0.253,clip]{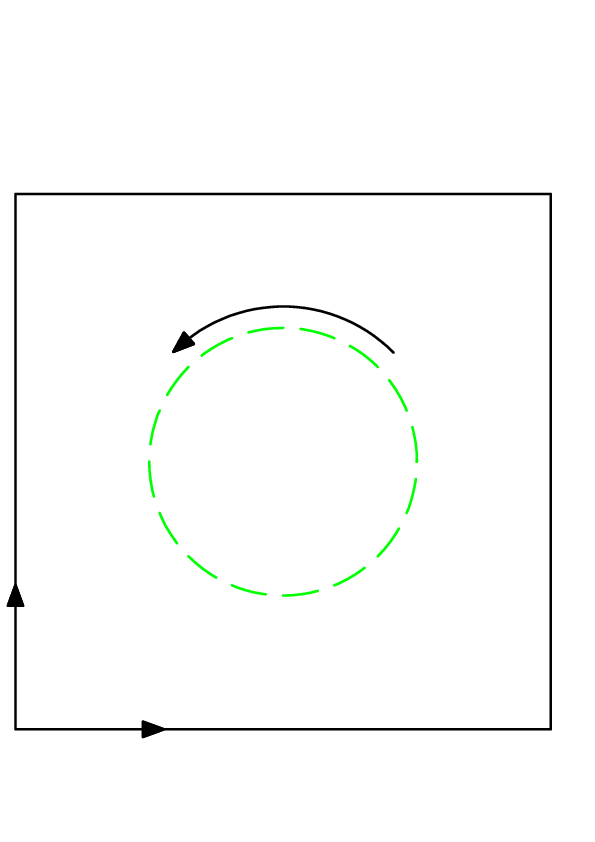}
		\begin{picture}(0,0)
		\put(-105,118){\small$|\tilde{\bm{u}}^f| = a(0.5 D)$}
		\put(-5,70){\small $\tilde{\bm{u}}^f_{ewns}$}
		\put(-5,55){\small $ = \bm{0.0}$}
		\put(-137,36){$y$}
		\put(-110,8){$x$}
		\end{picture}
		\caption{$\qquad \qquad$}
		\label{fig:schRot}
	\end{subfigure}
	\caption{Schematic diagram of convergence analysis: (a) Lid-driven cavity flow; (b) Rotating disk}
	\label{fig:errorScheme}
\end{figure}
\begin{figure} \centering
	\begin{subfigure}[b]{0.5\textwidth}
		\centering
		\hspace{-25pt}\includegraphics[trim=0.1cm 0.2cm 0.15cm 0.2cm,scale=0.22,clip]{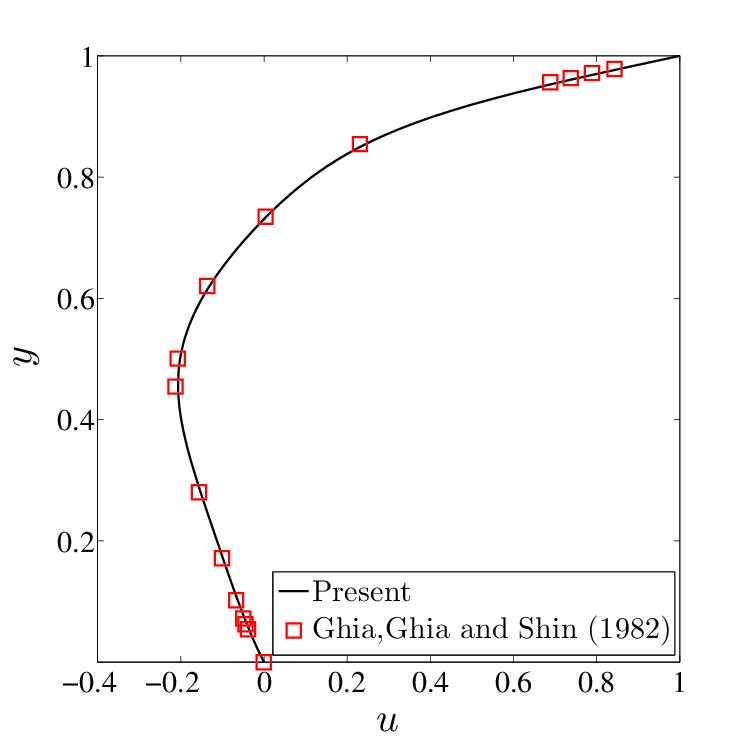}
		\caption{$\qquad$}
		\label{fig:u_vs_y}
	\end{subfigure}%
	\begin{subfigure}[b]{0.5\textwidth}	
		\centering
		\hspace{-25pt}\includegraphics[trim=0.1cm 0.1cm 0.1cm 0.1cm,scale=0.22,clip]{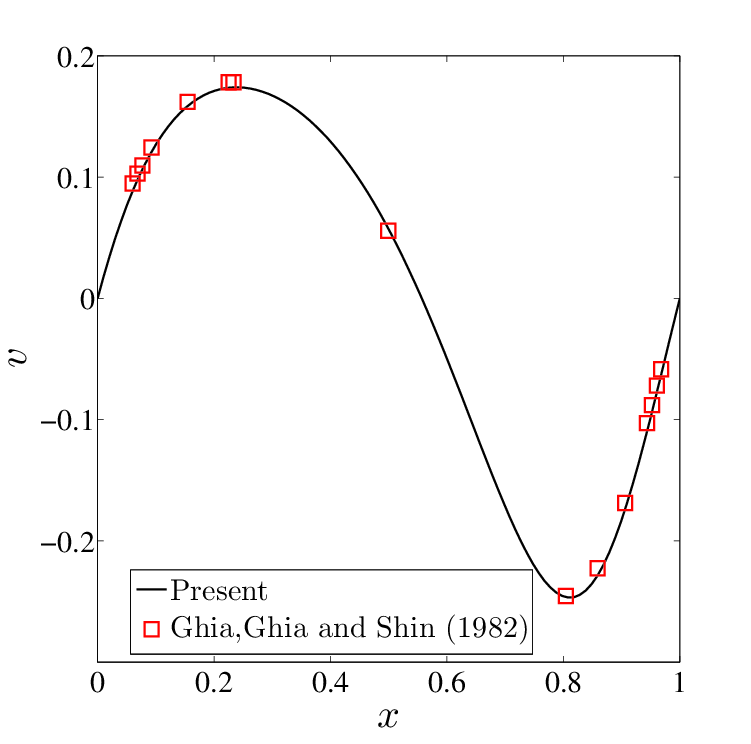}
		\caption{$\qquad$}
		\label{fig:v_vs_x}
	\end{subfigure}
	\caption{Velocity profile in classical lid-driven cavity flow at $Re=100$: (a) $u$ at $x=0.5$; (b) $v$ at $y=0.5$}
	\label{fig:validLid}
\end{figure}
\begin{figure} \centering
	\begin{subfigure}[b]{0.5\textwidth}	
		\centering
		\hspace{-25pt}\includegraphics[trim=0.1cm 0.1cm 0.1cm 0.1cm,scale=0.3,clip]{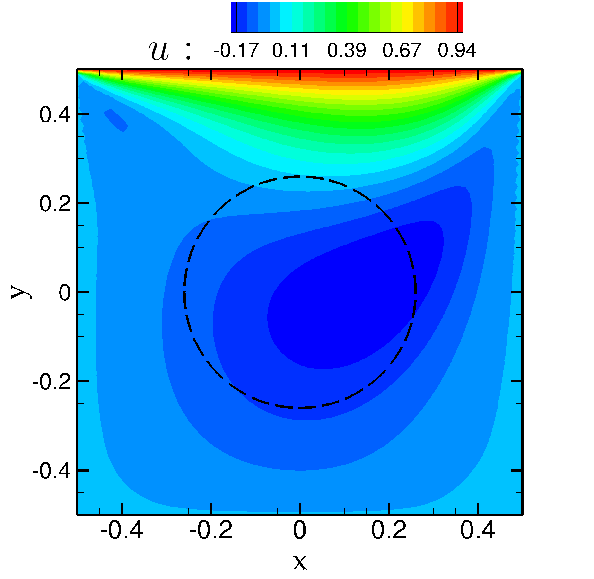}
		\caption{$\qquad$}
		\label{fig:uLid}
	\end{subfigure}%
	\begin{subfigure}[b]{0.5\textwidth}
		\centering
		\hspace{-25pt}\includegraphics[trim=0.1cm 0.1cm 0.1cm 0.1cm,scale=0.3,clip]{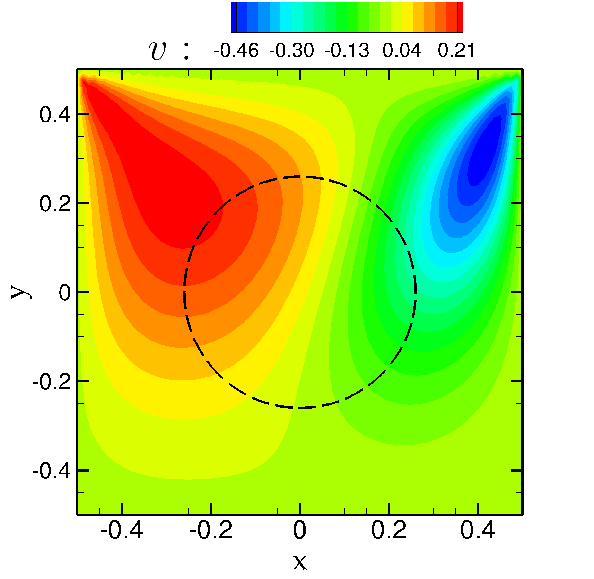}
		\caption{$\qquad$}
		\label{fig:vLid}
	\end{subfigure}
	\begin{subfigure}[b]{0.5\textwidth}	
		\centering
		\hspace{-25pt}\includegraphics[trim=0.1cm 0.1cm 0.1cm 0.1cm,scale=0.3,clip]{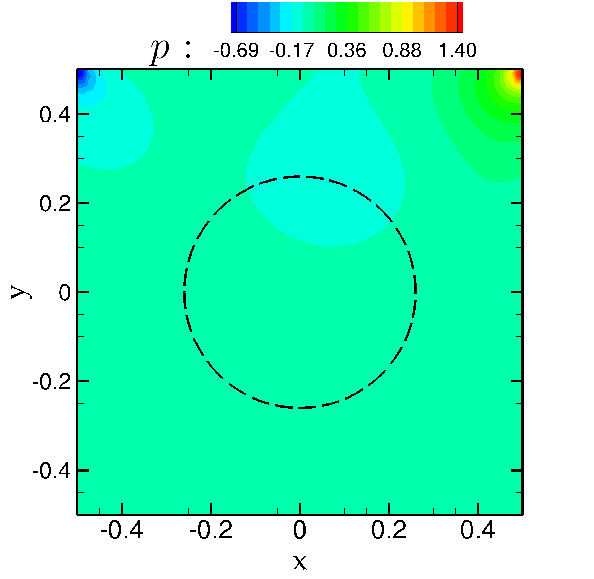}
		\caption{$\qquad$}
		\label{fig:pLid}
	\end{subfigure}%
	\begin{subfigure}[b]{0.5\textwidth}
		\centering
		\hspace{-25pt}\includegraphics[trim=0.1cm 0.1cm 0.1cm 0.1cm,scale=0.3,clip]{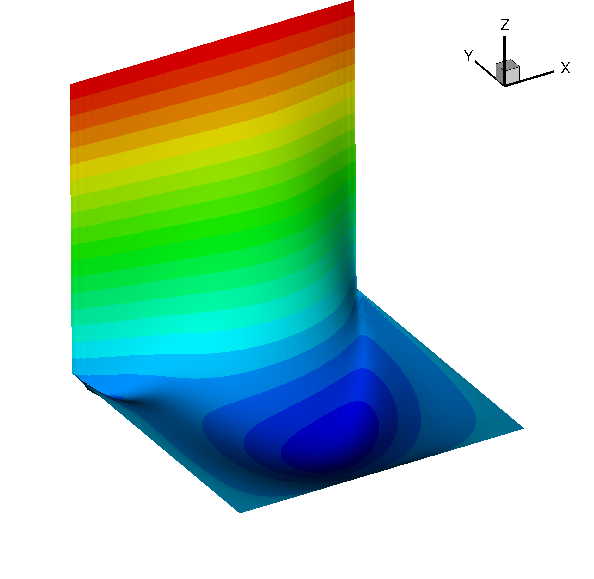}
		\caption{$\qquad$}
		\label{fig:isoLid}
	\end{subfigure}
	\caption{Contour plots for lid-driven cavity flow at $Re=100$: (a) x-component velocity field; (b) y-component velocity field; (c) pressure field; (d) 3D contour of u-component velocity field}
	\label{fig:conLid}
\end{figure}
\begin{figure} \centering
	\begin{subfigure}[b]{0.5\textwidth}	
		\centering
		\hspace{-25pt}\includegraphics[trim=0.1cm 0.1cm 0.1cm 0.1cm,scale=0.3,clip]{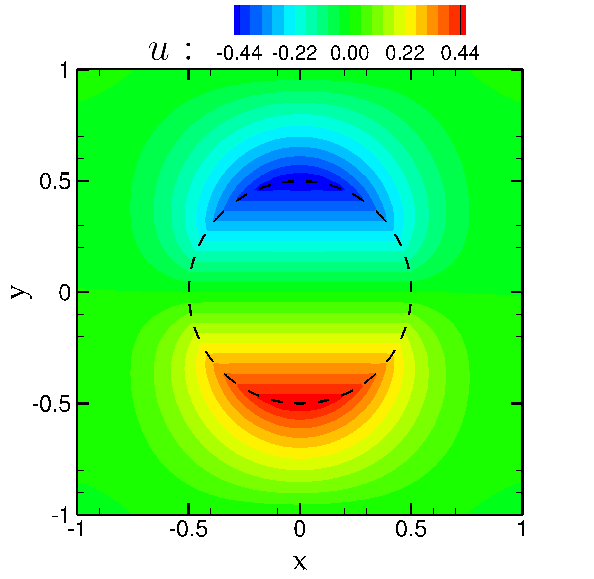}
		\caption{$\qquad$}
		\label{fig:uRot}
	\end{subfigure}%
	\begin{subfigure}[b]{0.5\textwidth}
		\centering
		\hspace{-25pt}\includegraphics[trim=0.1cm 0.1cm 0.1cm 0.1cm,scale=0.3,clip]{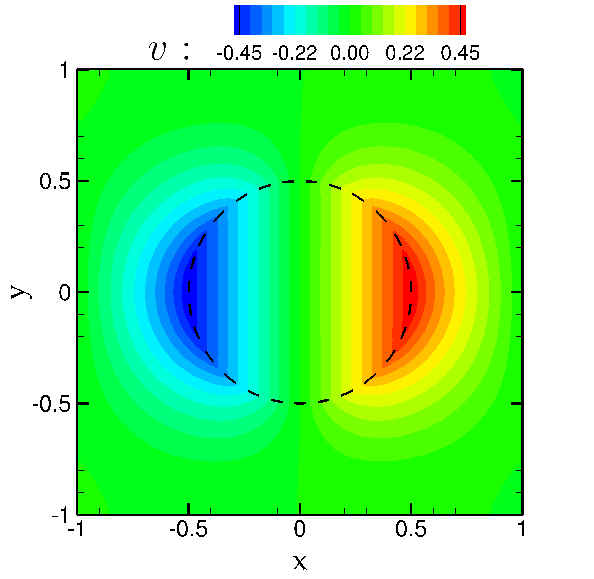}
		\caption{$\qquad$}
		\label{fig:vRot}
	\end{subfigure}
	\begin{subfigure}[b]{0.5\textwidth}	
		\centering
		\hspace{-25pt}\includegraphics[trim=0.1cm 0.1cm 0.1cm 0.1cm,scale=0.3,clip]{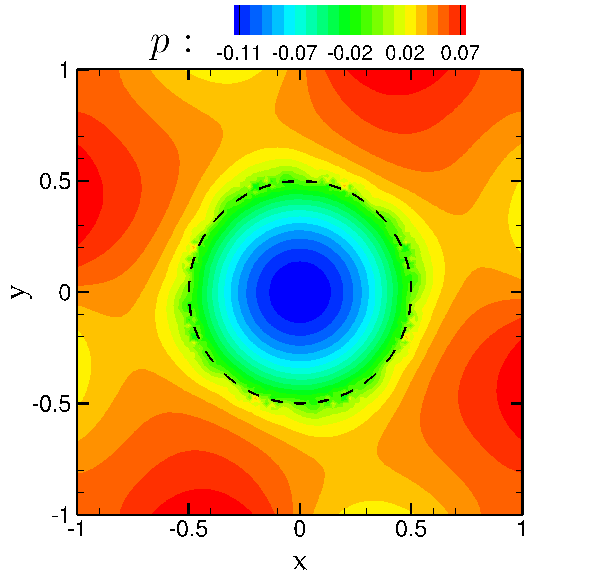}
		\caption{$\qquad$}
		\label{fig:pRot}
	\end{subfigure}%
	\begin{subfigure}[b]{0.5\textwidth}
		\centering
		\hspace{-25pt}\includegraphics[trim=0.1cm 0.1cm 0.11cm 0.1cm,scale=0.3,clip]{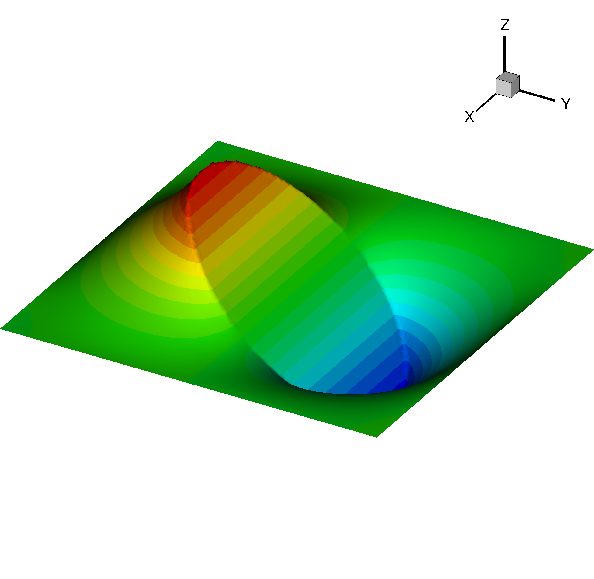}
		\caption{$\qquad$}
		\label{fig:isoRot}
	\end{subfigure}
	\caption{Contour plots for rotating disk at $Re=100$ and $a = 1.0$: (a) x-component velocity field; (b) y-component velocity field; (c) pressure field; (d) 3D contour of u-component velocity field}
	\label{fig:conRot}
\end{figure}

The convergence analyses of proposed PGQ are conducted via simulations of a classical lid-driven cavity flow and a rotating disk. The embedded interface is represented by a level-set function. In the lid-driven cavity flow, no Dirichlet boundary condition is weakly-imposed along the embedded interface, as shown in Fig.~\ref{fig:schLid}, where the subscripts $"e"$, $"w"$, $"n"$ and $"s"$ respectively refer to the east, west, north and south wall boundaries. Its objective is to get rid of influence by Nitsche's method and barely investigate the convergence rate of PGQ. The numerical results from the lid-driven cavity flow agree well with literature~\cite{Ghia1982Jocp}, as shown in Fig.~\ref{fig:validLid}. In Fig.~\ref{fig:conLid}, the resultant contours of velocity and pressure are smooth across the elements implemented with PGQ and Gauss quadrature numerical integrations. No odd values are observed in contours of variable across the embedded interface. It means that the proposed PGQ technique is well suited for working together with Gauss quadrature. On the other hand, a prescribed velocity is weakly-imposed in the rotating disk case Fig.~\ref{fig:schRot}. Similar to the lid-driven cavity flow, no odd value is observed in the contours of variable across the embedded interface in Fig.~\ref{fig:conRot}. In addition, prominent discontinuities in the value of pressure and the gradients of velocity are observed along the embedded interface, as shown in Fig.~\ref{fig:pRot} and~\ref{fig:isoRot} respectively. 
\begin{figure} \centering
	\begin{subfigure}[b]{0.5\textwidth}
		\centering
		\hspace{-25pt}\includegraphics[trim=0.1cm 0.2cm 0.15cm 0.2cm,scale=0.22,clip]{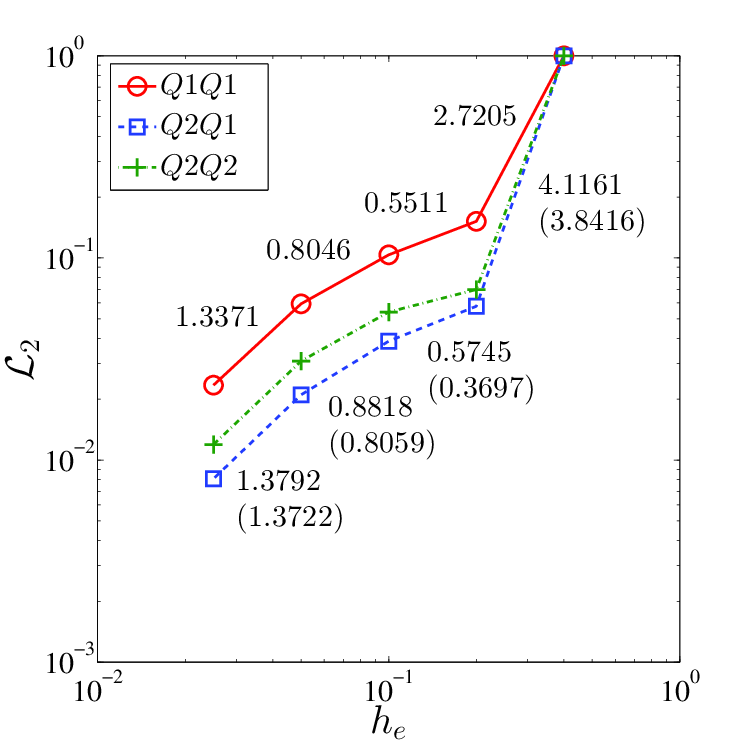}
		\caption{$\qquad$}
		\label{fig:eLD}
	\end{subfigure}%
	\begin{subfigure}[b]{0.5\textwidth}	
		\centering
		\hspace{-25pt}\includegraphics[trim=0.1cm 0.1cm 0.1cm 0.1cm,scale=0.22,clip]{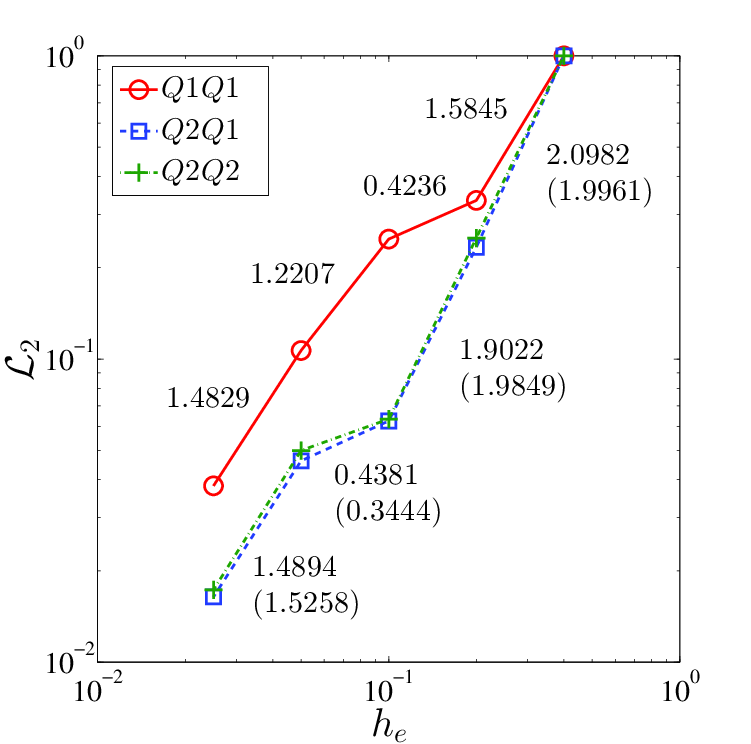}
		\caption{$\qquad$}
		\label{fig:eRotDisk}
	\end{subfigure}
\caption{Order of accuracy with respect to Eulerian grid refinement: (a) Lid-driven cavity flow at $Re=100$; (b) Rotating disk at $Re=20$ and $a = 1.0$}
	\label{fig:errorOrder}
\end{figure}

The convergence analyses of the lid-driven cavity flow and rotating disk are conducted for different element types, e.g., Q1Q1, Q2Q1 and Q2Q2. The results are plotted in Fig.~\ref{fig:errorOrder}, where $\mathcal{L}_2$ and $h_e$ denote the Euclidean 2 norm and the element length respectively. $\mathcal{L}_2$ norm is computed based on Eq.~\eqref{eq:l2}, in which $\bm{E}$ and $\varphi$ are the relative error vector and measured quantity respectively, e.g, x-component velocity. The superscript $"n"$ and subscript $"ref"$ respectively denote the number of background nodes along a side and the numerical results with a reference grid. The order of convergence is computed based on $\bm{E}$ in Eq.~\eqref{eq:order}. 
\begin{eqnarray}
||\bm{E}^n_{\varphi}||_{\mathcal{L}_2} &=& \sqrt{(\bm{E}^n_{\varphi})' \cdot \bm{E}^n_{\varphi}} \;\;;\qquad E^n_{\varphi} (i) = \varphi (i) - \varphi_{ref} (\bm{x} (i)) \label{eq:l2}\\
\text{order} &=& \frac{log(||\bm{E}^n_{\varphi}||_{\mathcal{L}_2}/||\bm{E}^{2n}_{\varphi}||_{\mathcal{L}_2})}{log(h^n_e/h^{2n}_e)} \label{eq:order}
\end{eqnarray}
The convergence rates are annotated in the plots. Those in parenthesis are associated with Q2Q2 element, the green line. In all simulations, the convergence rates for higher order element is approximately one-order higher than linear elements. In lid-driven cavity flow in Fig.~\ref{fig:eLD}, the convergence rates are approximately 2 and 2.8 for bi-linear and bi-quadratic elements respectively. By weakly-imposing a prescribed velocity along the embedded interface, the convergence rates for bi-linear and bi-quadratic elements are approximately 1.45 and 2.0 respectively. It is approximately half-order lower than the case of lid-driven cavity flow. 

\section{Numerical examples and Validations} \label{sec:NE}
In this section, a number of representative numerical examples are presented to assess the accuracy and robustness of the proposed PGQ technique. The performed simulations are (a) a stationary cylinder in cross-flow, (b) a rotating cylinder in cross-flow, (c) a freely-vibrating cylinder in cross-flow, (d) a free-falling particle and (e) six free-falling particles.  

\subsection{Stationary cylinder in cross-flow}
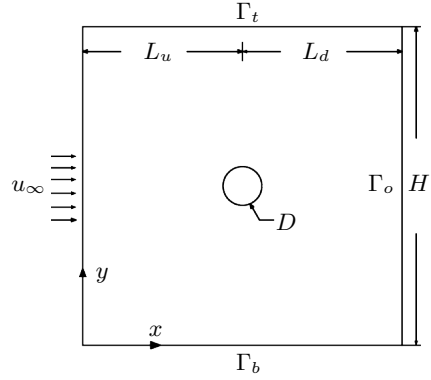
\begin{figure} \centering
	\hspace{-25pt}\includegraphics[trim=0.1cm 4cm 1cm 6.5cm,scale=0.25,clip]{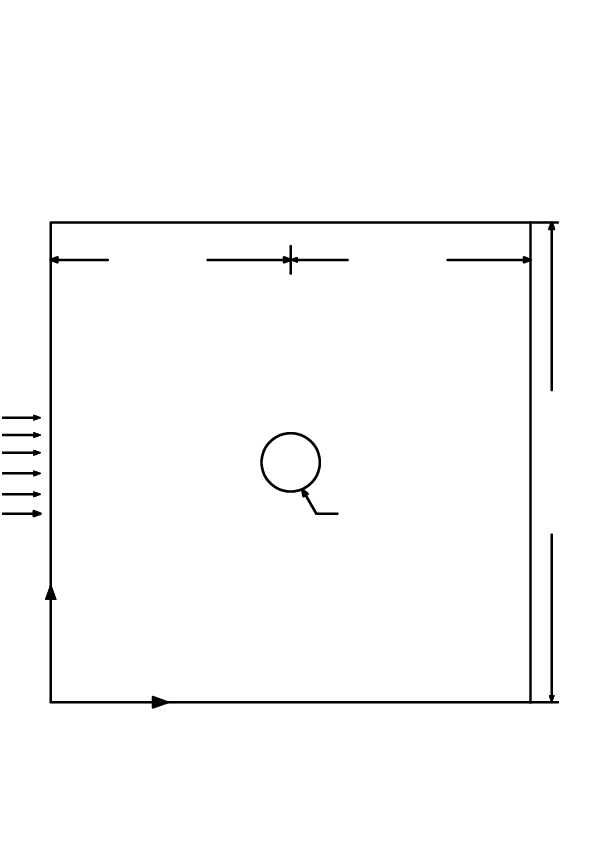}
	\begin{picture}(0,0)
	\put(-108,10){\small $x$}
	\put(-128,32){\small $y$}
	\put(-160,65){\small $u_{\infty}$}
	\put(-60,50){\small $D$}
	\put(-110,115){\small $L_u$}
	\put(-50,115){\small $L_d$}
	\put(-10,65){\small $H$}
	\put(-25,65){\small $\Gamma_o$}
	\put(-75,130){\small $\Gamma_t$}
	\put(-75,-3){\small $\Gamma_b$}
	\end{picture}
	\caption{Schematic diagram of a stationary cylinder in cross-flow}
	\label{fig:fixCyl}
\end{figure}
\begin{table*}
	\setlength{\tabcolsep}{12pt}
	\renewcommand{\arraystretch}{2}
	\begin{center}
		\def\arraystretch{0.8}
		\begin{tabular}{l p{1cm} p{1cm} p{1cm}}
			\hline
			&  & $L/D$ &  $C_d$\\
			\hline
			\multirow{4}{*}{$Re=20$}  & \multicolumn{1}{l}{Tritton \cite{tritton1959experiments}} & \multicolumn{1}{l}{---} & \multicolumn{1}{l}{2.22}\\
			& \multicolumn{1}{l}{Coutanceau and Bouard \cite{coutanceau1977experimental}} &  \multicolumn{1}{l}{0.73} & \multicolumn{1}{l}{---} \\
			& \multicolumn{1}{l}{Calhoun \cite{calhoun2002cartesian}} &  \multicolumn{1}{l}{0.91} & \multicolumn{1}{l}{2.19} \\
			& \multicolumn{1}{l}{Russell and Wang \cite{russell2003cartesian}} &  \multicolumn{1}{l}{0.94} & \multicolumn{1}{l}{2.13} \\
			& \multicolumn{1}{l}{Li et al. \cite{li2016immersed}} &  \multicolumn{1}{l}{0.931} & \multicolumn{1}{l}{2.062} \\
			& \parbox[t]{5cm}{Present} & \multicolumn{1}{l}{0.94} &\multicolumn{1}{l}{2.171} \\
			\cline{1-4}
			\multirow{4}{*}{$Re=40$}  & \multicolumn{1}{l}{Tritton \cite{tritton1959experiments}} & \multicolumn{1}{l}{---} & \multicolumn{1}{l}{1.48}\\
			& \multicolumn{1}{l}{Coutanceau and Bouard \cite{coutanceau1977experimental}} &  \multicolumn{1}{l}{1.89} & \multicolumn{1}{l}{---} \\
			& \multicolumn{1}{l}{Calhoun \cite{calhoun2002cartesian}} &  \multicolumn{1}{l}{2.18} & \multicolumn{1}{l}{1.62} \\
			& \multicolumn{1}{l}{Russell and Wang \cite{russell2003cartesian}} &  \multicolumn{1}{l}{2.29} & \multicolumn{1}{l}{1.60} \\
			& \multicolumn{1}{l}{Li et al. \cite{li2016immersed}} &  \multicolumn{1}{l}{2.24} & \multicolumn{1}{l}{1.569} \\
			& \parbox[t]{5cm}{Present} & \multicolumn{1}{l}{2.27} &\multicolumn{1}{l}{1.608} \\
			\hline
		\end{tabular}
	\end{center}
	\caption{Flow around a fixed circular cylinder: $L/D$ and $C_d$ for $Re=20$ and $40$}
	\label{tab:statVal1}
\end{table*}
\begin{table*}
\setlength{\tabcolsep}{12pt}
\renewcommand{\arraystretch}{2}
\begin{center}
\def\arraystretch{0.8}
\begin{tabular}{l p{1cm} p{1cm} p{1cm} p{1cm}}
\hline
&  & $C^{mean}_d$ &  $C^{rms}_l$ & $S_t$\\
\hline
\multirow{4}{*}{$Re=100$}  & \multicolumn{1}{l}{Braza et al. \cite{braza1986numerical}} & \multicolumn{1}{l}{1.364} & \multicolumn{1}{l}{$\pm$0.25} & \multicolumn{1}{l}{---}\\
& \multicolumn{1}{l}{Liu et al. \cite{liu1998preconditioned}} &  \multicolumn{1}{l}{1.350} & \multicolumn{1}{l}{$\pm$0.339} & \multicolumn{1}{l}{0.164}  \\
& \multicolumn{1}{l}{Calhoun \cite{calhoun2002cartesian}} &  \multicolumn{1}{l}{1.330} & \multicolumn{1}{l}{$\pm$0.298} & \multicolumn{1}{l}{0.175}  \\
& \multicolumn{1}{l}{Russell and Wang \cite{russell2003cartesian}} &  \multicolumn{1}{l}{1.380} & \multicolumn{1}{l}{$\pm$0.300} & \multicolumn{1}{l}{0.169}  \\
& \multicolumn{1}{l}{Li et al. \cite{li2016immersed}} &  \multicolumn{1}{l}{1.301} & \multicolumn{1}{l}{$\pm$0.324} & \multicolumn{1}{l}{0.167}  \\
& \multicolumn{1}{l}{Kadapa et al. \cite{kadapa2016fictitious}} &  \multicolumn{1}{l}{1.390} & \multicolumn{1}{l}{$\pm$0.339} & \multicolumn{1}{l}{0.166}  \\
& \parbox[t]{5cm}{Present} & \multicolumn{1}{l}{1.365} &\multicolumn{1}{l}{$\pm$0.301} & \multicolumn{1}{l}{0.164}  \\
\cline{1-5}
\multirow{4}{*}{$Re=200$}  & \multicolumn{1}{l}{Braza et al. \cite{braza1986numerical}} & \multicolumn{1}{l}{1.40} & \multicolumn{1}{l}{$\pm$0.75} & \multicolumn{1}{l}{---}\\
& \multicolumn{1}{l}{Liu et al. \cite{liu1998preconditioned}} &  \multicolumn{1}{l}{1.310} & \multicolumn{1}{l}{$\pm$0.69} & \multicolumn{1}{l}{0.192}  \\
& \multicolumn{1}{l}{Calhoun \cite{calhoun2002cartesian}} &  \multicolumn{1}{l}{1.172} & \multicolumn{1}{l}{$\pm$0.594} & \multicolumn{1}{l}{0.202}  \\
& \multicolumn{1}{l}{Russell and Wang \cite{russell2003cartesian}} &  \multicolumn{1}{l}{1.390} & \multicolumn{1}{l}{$\pm$0.50} & \multicolumn{1}{l}{0.195}  \\
& \multicolumn{1}{l}{Li et al. \cite{li2016immersed}} &  \multicolumn{1}{l}{1.307} & \multicolumn{1}{l}{$\pm$0.419} & \multicolumn{1}{l}{0.192}  \\
& \multicolumn{1}{l}{Kadapa et al. \cite{kadapa2016fictitious}} &  \multicolumn{1}{l}{1.42} & \multicolumn{1}{l}{$\pm$0.594} & \multicolumn{1}{l}{0.194}  \\
& \parbox[t]{5cm}{Present} & \multicolumn{1}{l}{1.372} &\multicolumn{1}{l}{$\pm$0.648} & \multicolumn{1}{l}{0.194}  \\
\hline
\end{tabular}
\end{center}
\caption{Flow around a fixed circular cylinder: $C^{mean}_d$, $C^{rms}_l$ and $St$ for $Re=100$ and $200$}
\label{tab:statVal2}
\end{table*}
\begin{figure} \centering
	\begin{subfigure}[b]{0.5\textwidth}	
		\centering
		\hspace{-25pt}\includegraphics[trim=0.1cm 0.1cm 0.1cm 0.1cm,scale=0.3,clip]{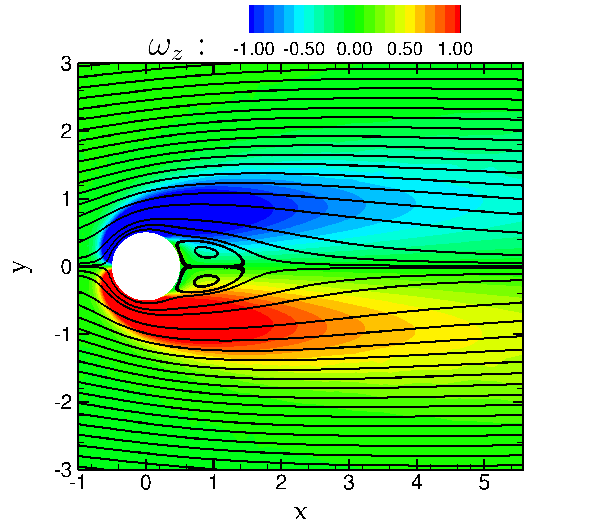}
		\caption{$\qquad$}
		\label{fig:re20wz}
	\end{subfigure}%
	\begin{subfigure}[b]{0.5\textwidth}
		\centering
		\hspace{-25pt}\includegraphics[trim=0.1cm 0.1cm 0.1cm 0.1cm,scale=0.3,clip]{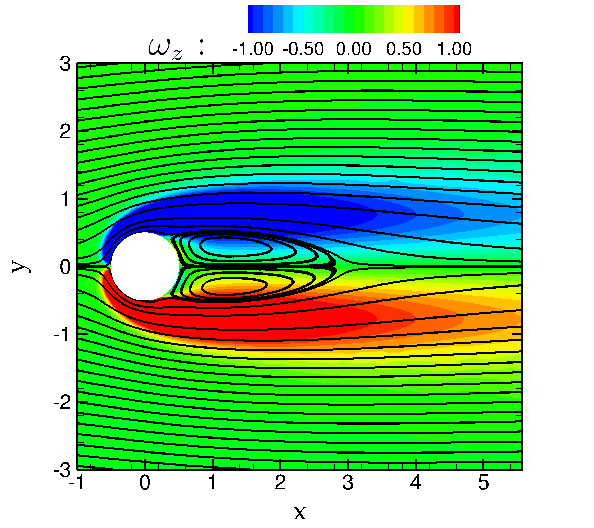}
		\caption{$\qquad$}
		\label{fig:re40wz}
	\end{subfigure}
	\caption{$\omega_z$ contour and streamline plot of a fixed circular cylinder: (a) $Re=20$; (b) $Re=40$}
	\label{fig:re2040wz}
\end{figure}
\begin{figure} \centering
\begin{subfigure}[b]{0.5\textwidth}	
\centering
\hspace{-25pt}\includegraphics[trim=0.1cm 0.1cm 0.1cm 0.1cm,scale=0.3,clip]{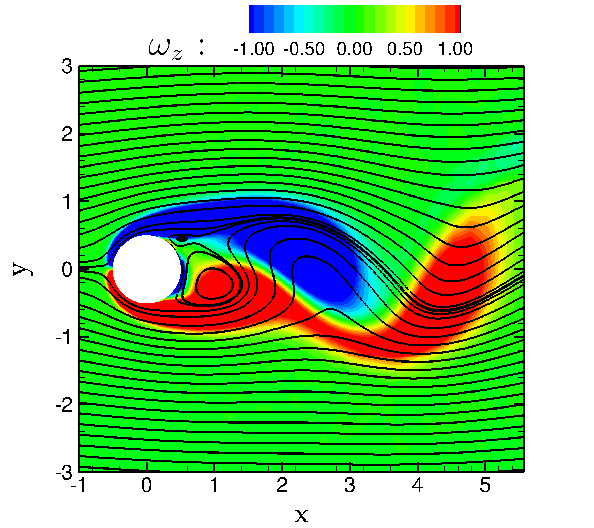}
\caption{$\qquad$}
\label{fig:re100wz}
\end{subfigure}%
\begin{subfigure}[b]{0.5\textwidth}
\centering
\hspace{-25pt}\includegraphics[trim=0.1cm 0.1cm 0.1cm 0.1cm,scale=0.3,clip]{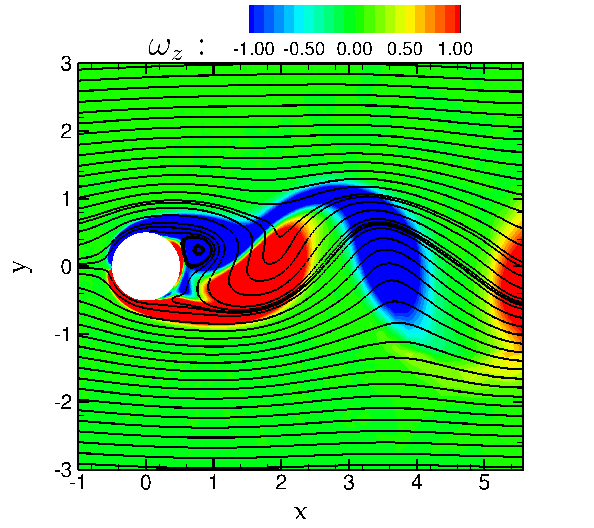}
\caption{$\qquad$}
\label{fig:re200wz}
\end{subfigure}
\caption{$\omega_z$ contour and streamline plot of a fixed circular cylinder: (a) $Re=100$; (b) $Re=200$}
\label{fig:re100200wz}
\end{figure}
The flow around a stationary cylinder in laminar flow, $Re \leq 200$, is a classical benchmark example. Its schematic diagram is shown in Fig.~\ref{fig:fixCyl}, where $u_{\infty} = 1.0$, $D=1.0$, $L_u=50D$, $L_d=50D$ and $H=100D$ denote the free stream velocity, diameter of cylinder, upstream length, downstream length and width of fluid domain. Traction free boundary condition is imposed on domain boundaries $\Gamma_o$, $\Gamma_t$ and $\Gamma_b$. The fluid density $\rho^f = 1.0$ and dynamic viscosity $\mu = 0.01$ is chosen for simulation. 

The numerical results are compared with literature and summarized in Tab.~\ref{tab:statVal1} and~\ref{tab:statVal2}. It shows the numerical results obtained from PGQ agree well with literature. The corresponding contour of z-component vorticity $\omega_z$ are plotted in Fig.~\ref{fig:re2040wz} and~\ref{fig:re100200wz}.

\subsection{Rotating cylinder in cross-flow}
\begin{figure} \centering
	\begin{subfigure}[b]{0.5\textwidth}	
		\centering
		\hspace{-25pt}\includegraphics[trim=0.1cm 3cm 1cm 5.5cm,scale=0.25,clip]{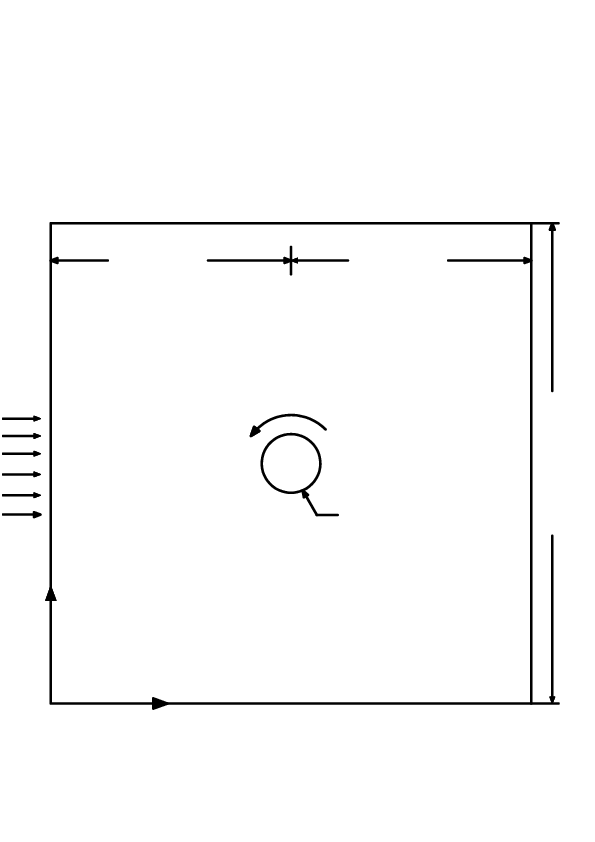}
		\begin{picture}(0,0)
		\put(-108,18){\small $x$}
		\put(-128,38){\small $y$}
		\put(-160,72){\small $u_{\infty}$}
		\put(-60,57){\small $D$}
		\put(-110,121){\small $L_u$}
		\put(-50,121){\small $L_d$}
		\put(-10,72){\small $H$}
		\put(-25,72){\small $\Gamma_o$}
		\put(-75,137){\small $\Gamma_t$}
		\put(-75,2){\small $\Gamma_b$}
		\put(-85,91){\footnotesize $a (0.5D)$}
		\end{picture}
		\caption{$\qquad$}
		\label{fig:cyl2}
	\end{subfigure}%
	\begin{subfigure}[b]{0.5\textwidth}
		\centering
		\hspace{-25pt}\includegraphics[trim=0.1cm 0.1cm 0.1cm 0.1cm,scale=0.3,clip]{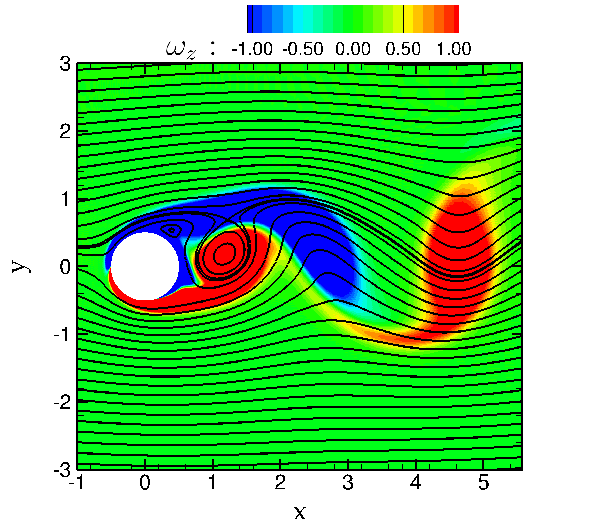}
		\caption{$\qquad$}
		\label{fig:re200a1con}
	\end{subfigure}
	\begin{subfigure}[b]{1.0\textwidth}	
		\centering
		\hspace{-25pt}\includegraphics[trim=0.01cm 0.1cm 0.1cm 0.1cm,scale=0.22,clip]{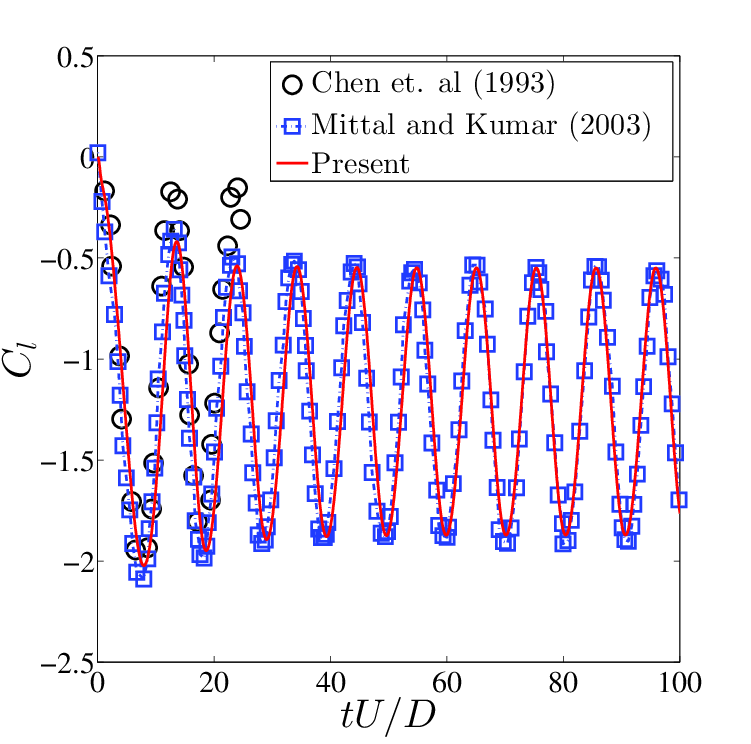}
		\caption{$\qquad$}
		\label{fig:re200a1Cl}
	\end{subfigure}
	\caption{Rotating cylinder in cross-flow at $Re=200$ and $a=1.0$: (a) Schematic diagram; (b) $\omega_z$ contour and streamline plot; (c) Time trace of lift coefficient}
	\label{fig:re200a1}
\end{figure}
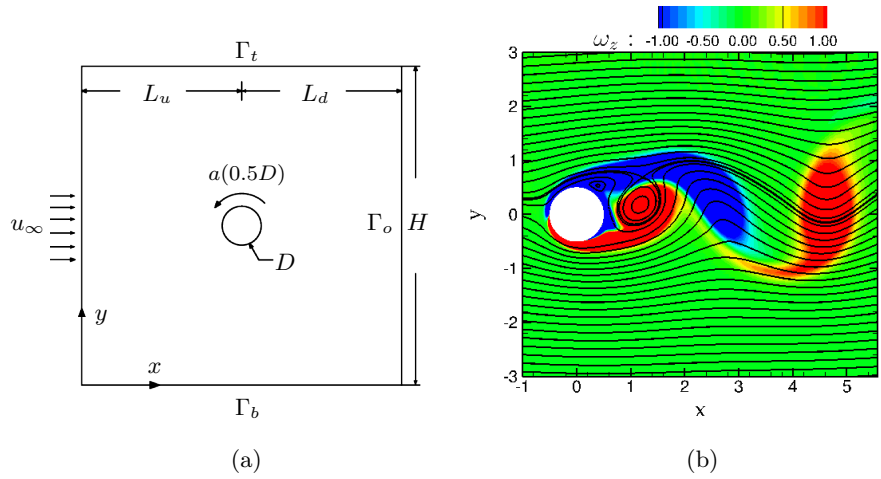
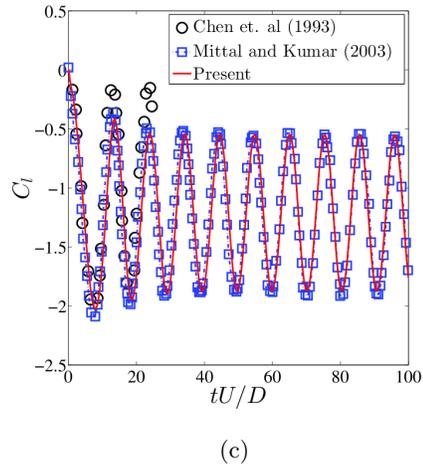
To simulate a rotating cylinder, a prescribed velocity $\tilde{\bm{u}}$ is imposed along the embedded interface. Its schematic diagram is shown in Fig.~\ref{fig:cyl2}. $\tilde{\bm{u}}$ is computed as $[a (0.5D)] \bm{n}$, where $a=1.0$, $D=1.0$ and $\bm{n}$ respectively are angular velocity, diameter of cylinder and basis vector. The corresponding contour of $\omega_z$ is plotted in Fig.~\ref{fig:re200a1con}. The response of lift coefficient agrees well with results from literature~\cite{Chen1993JoFM,Mittal2003JoFM}, as shown in Fig.~\ref{fig:re200a1Cl}. 

The impulsive initial data poses a challenge of convergence in the initial stage of simulation. To obtain a good convergence rate, the field data of a stationary cylinder is chosen as the initial condition. Since a second-order Generalized-$\alpha$ temporal integration scheme is implemented, accurate numerical results can be obtained at a relative larger time step, e.g., $dt = 0.02$

\subsection{Vibrating cylinder in cross-flow}
\begin{figure} \centering
	\begin{subfigure}[b]{0.5\textwidth}
		\centering
		\hspace{-25pt}\includegraphics[trim=0.1cm 4cm 1cm 5.5cm,scale=0.25,clip]{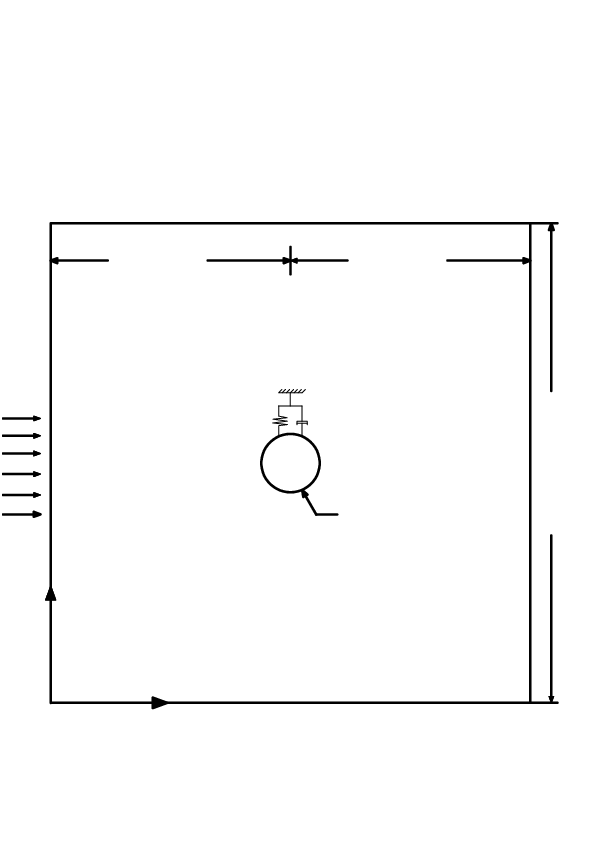}
		\begin{picture}(0,0)
		\put(-108,10){\small $x$}
		\put(-128,32){\small $y$}
		\put(-160,65){\small $u_{\infty}$}
		\put(-60,50){\small $D$}
		\put(-110,115){\small $L_u$}
		\put(-50,115){\small $L_d$}
		\put(-10,65){\small $H$}
		\put(-25,65){\small $\Gamma_o$}
		\put(-75,130){\small $\Gamma_t$}
		\put(-75,-3){\small $\Gamma_b$}
		\end{picture}
		\caption{$\qquad$}
		\label{fig:cyl3}
	\end{subfigure}%
	\begin{subfigure}[b]{0.5\textwidth}
		\centering
		\hspace{-25pt}\includegraphics[trim=0.1cm 4cm 1cm 5.5cm,scale=0.25,clip]{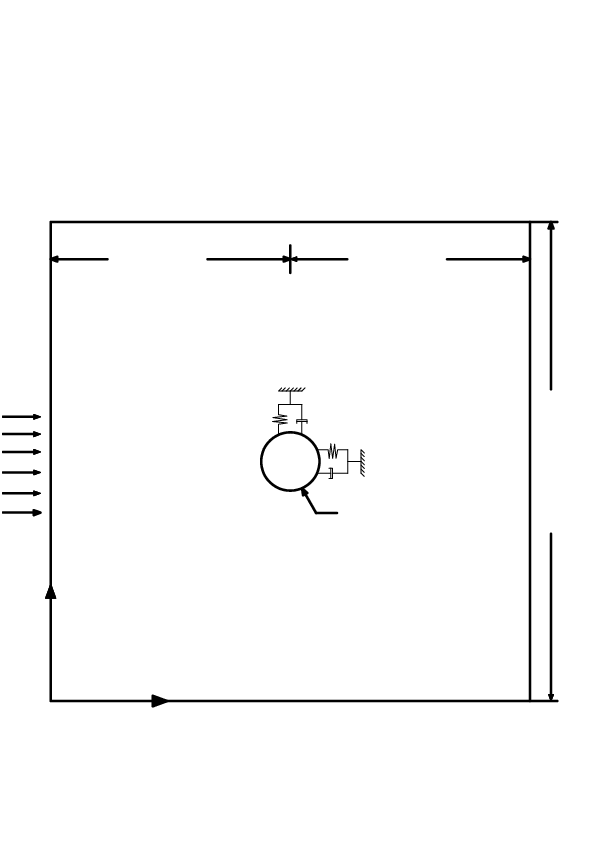}
		\begin{picture}(0,0)
		\put(-108,10){\small $x$}
		\put(-128,32){\small $y$}
		\put(-160,65){\small $u_{\infty}$}
		\put(-60,50){\small $D$}
		\put(-110,115){\small $L_u$}
		\put(-50,115){\small $L_d$}
		\put(-10,65){\small $H$}
		\put(-25,65){\small $\Gamma_o$}
		\put(-75,130){\small $\Gamma_t$}
		\put(-75,-3){\small $\Gamma_b$}
		\end{picture}
		\caption{$\qquad$}
		\label{fig:cyl4}
	\end{subfigure}
	\caption{Schematic diagrams of vibrating cylinder in cross-flow: (a) a transverse-vibrating (1-DoFs) cylinder; (d) a freely-vibrating (2-DoFs) cylinder}
	\label{fig:vCyl}
\end{figure}
\begin{figure} \centering
	\begin{subfigure}[b]{0.5\textwidth}	
		\centering
		\hspace{-25pt}\includegraphics[trim=0.1cm 0.1cm 0.1cm 0.1cm,scale=0.22,clip]{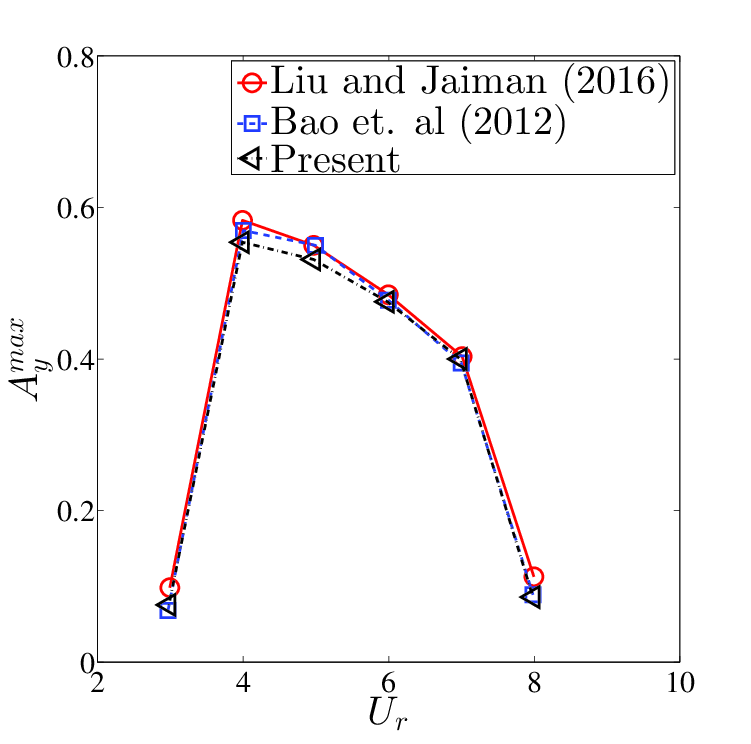}
		\caption{$\qquad$}
		\label{fig:Ay}
	\end{subfigure}%
	\begin{subfigure}[b]{0.5\textwidth}
		\centering
		\hspace{-25pt}\includegraphics[trim=0.1cm 0.1cm 0.1cm 0.1cm,scale=0.22,clip]{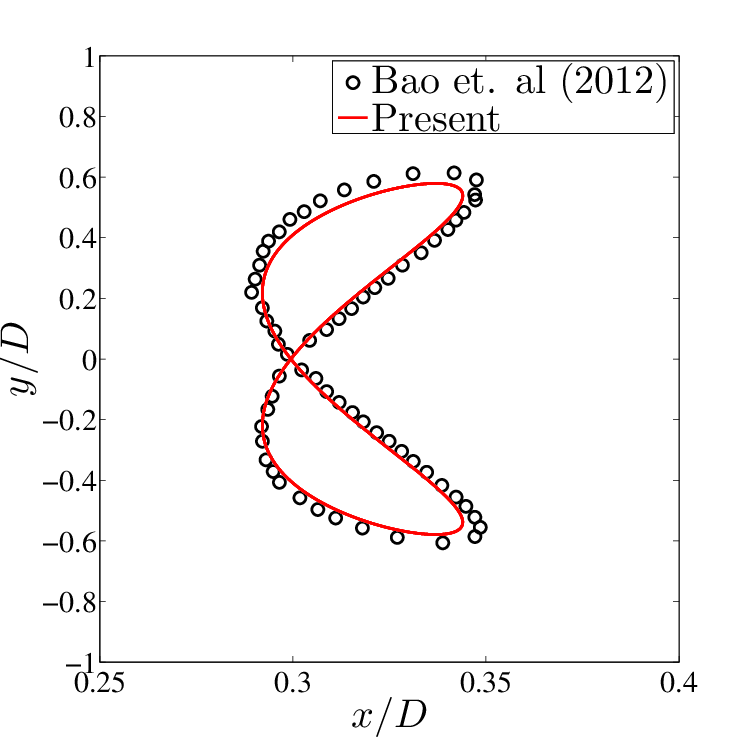}
		\caption{$\qquad$}
		\label{fig:traj}
	\end{subfigure}
	\begin{subfigure}[b]{0.5\textwidth}	
		\centering
		\hspace{-25pt}\includegraphics[trim=0.1cm 0.1cm 0.1cm 0.1cm,scale=0.3,clip]{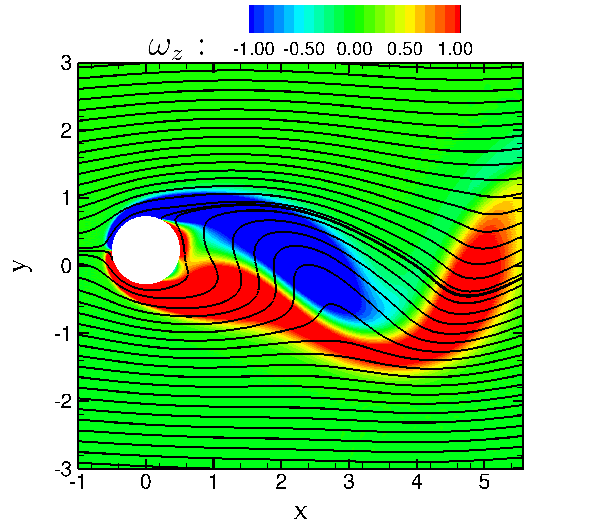}
		\caption{$\qquad$}
		\label{fig:wzRe100r7}
	\end{subfigure}%
	\begin{subfigure}[b]{0.5\textwidth}
		\centering
		\hspace{-25pt}\includegraphics[trim=0.1cm 0.1cm 0.1cm 0.1cm,scale=0.3,clip]{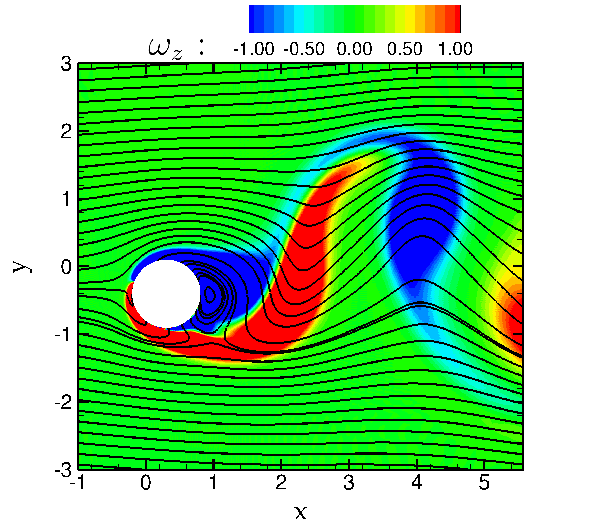}
		\caption{$\qquad$}
		\label{fig:wzRe150r5}
	\end{subfigure}
	\caption{Vibrating cylinder in cross-flow: (a,c) $Re=100$, $m^*=10.0$, $\zeta=0.01$, $U_r = 7.0$ and vibrating in y direction; (b,d) $Re=150$, $m^*=2.55$, $\zeta=0.0$, $U_r=5.0$ and vibrating in x and y directions}
	\label{fig:viv}
\end{figure}
In this section, two types of vibrating cylinder is chosen as benchmark examples, e.g., transversely-vibrating (1-DoFs) cylinder Fig.~\ref{fig:cyl3} and freely-vibrating (2-DoFs) cylinder in x and y directions Fig.~\ref{fig:cyl4}. For transversely-vibrating cylinder cases, $Re=100$, $m^*=10.0$, $\zeta = 0.01$ and $U_r \in$ [3,8] are chosen to set up the simulations. The obtained numerical results in Fig.~\ref{fig:Ay} show a good agreement with literature~\cite{liu2016interaction,Bao2012JoFaS}. In freely-vibrating cylinder case, the cylinder can vibrate in both x and y directions. A representative case is chosen for validation at $Re = 150$, $m^* = 2.55$, $\zeta=0.0$ and $U_r = 5.0$. Its trajectory results in Fig.~\ref{fig:traj} match well with literature~\cite{Bao2012JoFaS}. The contours of $\omega_z$ for representative cases are plotted in Fig.~\ref{fig:wzRe100r7} and~\ref{fig:wzRe150r5} respectively.

\subsection{Free-falling: a single particle} \label{sec:1fall}
\begin{figure} \centering
	\begin{subfigure}[b]{0.5\textwidth}	
		\centering
		\hspace{-25pt}\includegraphics[trim=0.1cm 0.1cm 0.1cm 0.1cm,scale=0.28,clip]{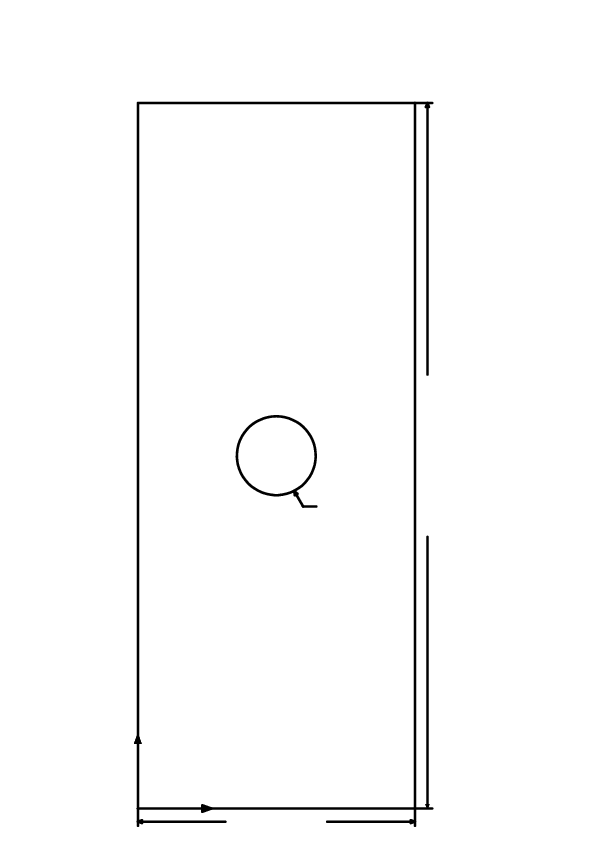}
		\begin{picture}(0,0)
		\put(-80,90){\small $D$}
		\put(-95,1){\small $2$}
		\put(-50,103){\small $6$}
		\put(-115,12){\small $x$}
		\put(-128,27){\small $y$}
		\put(-93,90){\vector(0,-1){20}}
		\put(-109,60){\footnotesize $g=9.81$}
		\put(-160,110){\small $\tilde{\bm{u}}^f_{ews}$}
		\put(-160,100){\small $=\bm{0.0}$}
		\put(-110,210){\small $\tilde{\bm{h}}^f = \bm{0.0}$}
		\end{picture}
		\caption{$\qquad$}
		\label{fig:freeschm}
	\end{subfigure}%
	\begin{subfigure}[b]{0.5\textwidth}	
		\centering
		\hspace{-25pt}\includegraphics[trim=0.1cm 2.55cm 0.1cm 2cm,scale=0.34,clip]{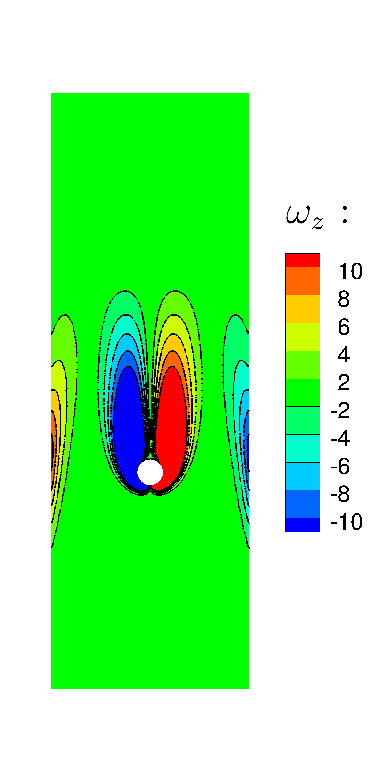}
		\caption{$\qquad$}
		\label{fig:wzFree}
	\end{subfigure}
	\begin{subfigure}[b]{0.5\textwidth}	
		\centering
		\hspace{-25pt}\includegraphics[trim=0.1cm 0.1cm 0.1cm 0.1cm,scale=0.22,clip]{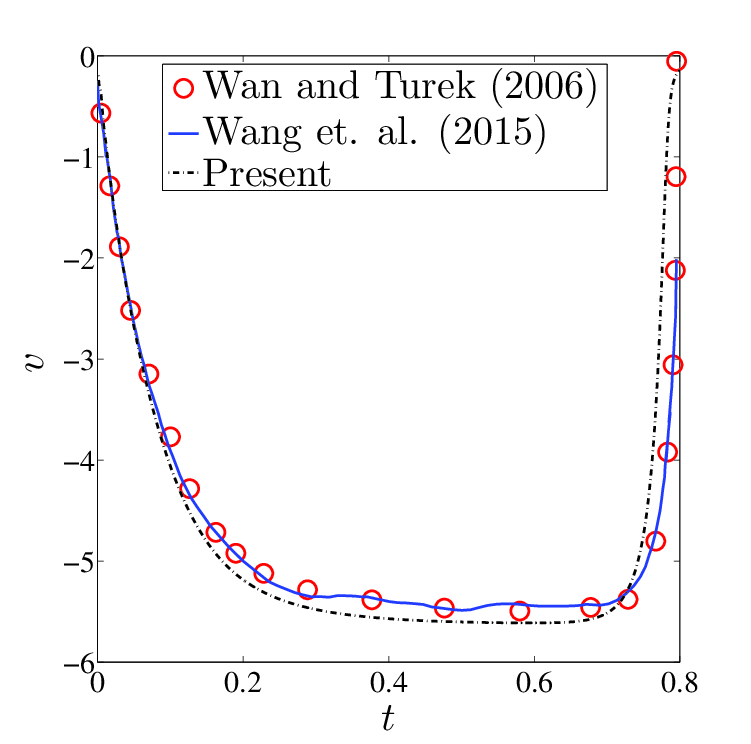}
		\caption{$\qquad$}
		\label{fig:vFree}
	\end{subfigure}%
	\begin{subfigure}[b]{0.5\textwidth}
		\centering
		\hspace{-25pt}\includegraphics[trim=0.1cm 0.1cm 0.1cm 0.1cm,scale=0.22,clip]{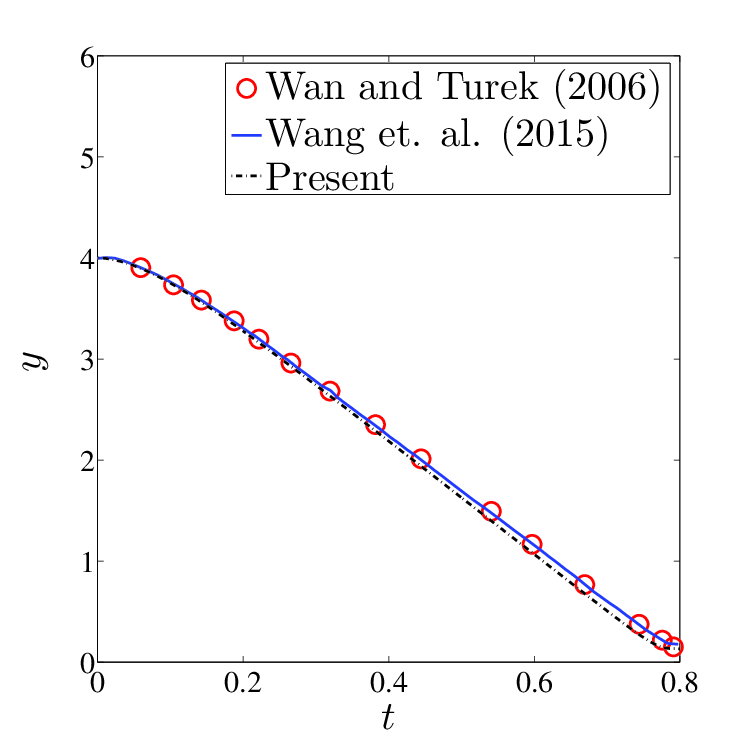}
		\caption{$\qquad$}
		\label{fig:yFree}
	\end{subfigure}
	\caption{Free-falling particle at $Re=13.75$, $m^*=1.25$ and $\zeta=0.01$: (a) schematic diagram; (b) $\omega_z$ contour plot at $t=0.4$; (c) time trace of y-component velocity; (d) time trace of y-component displacement}
	\label{fig:freeFall}
\end{figure}
Sedimentation is a classical benchmark example for fictitious domain methods. In this example, a circular particle is free-falling under gravitational force in an incompressible Newtonian fluid. The particle is accelerated at rest and subsequently achieve a terminal velocity $\bm{u}^s_t$. The chosen parameters in the simulation are $m^*=1.25$, $\zeta=0.01$, $\rho^f = 1000$, $\mu = 0.01$ and $D=0.25$. 

The schematic diagram is shown in Fig.~\ref{fig:freeschm}. The subscript $"e"$, $"w"$ and $"s"$ denotes the east, west and south wall boundary respectively. "no-slip" boundary condition is imposed on the east, west and south walls $\tilde{\bm{u}}^f_{ews} = \bm{0.0}$. Traction free boundary condition is imposed on the output as $\tilde{\bm{h}}^f = \bm{0.0}$. The particle falls from the rest at $[x,y] = [1,4]$. The contour of $\omega_z$ is plotted in Fig.~\ref{fig:wzFree}. The numerical results is compared with literature~\cite{Wan2006IJfNMiF,wang2015immersed} in Fig.~\ref{fig:vFree} and~\ref{fig:yFree}. The obtained numerical results can match with literature well.

\subsection{Free falling: 6 particles} \label{sec:2fall}
\begin{figure} \centering
	\begin{subfigure}[b]{0.5\textwidth}	
		\centering
		\hspace{-25pt}\includegraphics[trim=3cm 3cm 3cm 0.1cm,scale=0.3,clip]{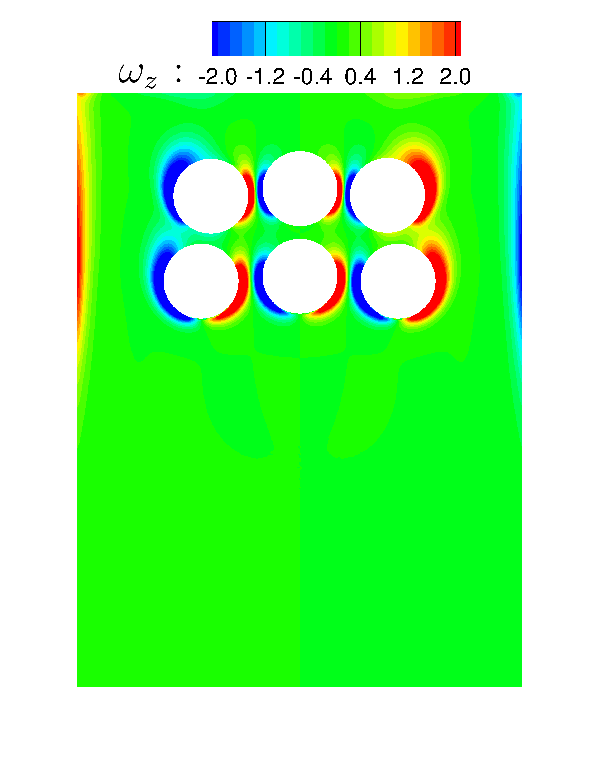}
		\caption{$t=3$}
		\label{fig:6fall1}
	\end{subfigure}%
	\begin{subfigure}[b]{0.5\textwidth}	
		\centering
		\hspace{-25pt}\includegraphics[trim=3cm 3cm 3cm 0.1cm,scale=0.3,clip]{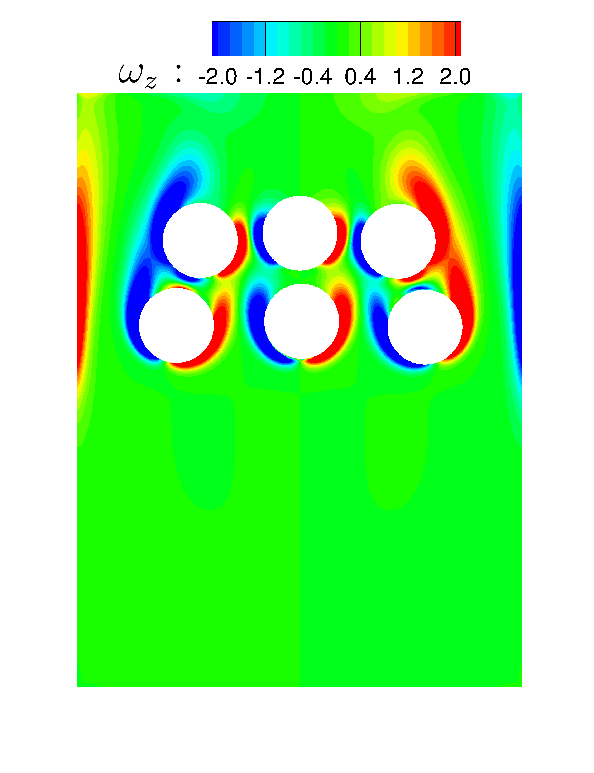}
		\caption{$t=6$}
		\label{fig:6fall2}
	\end{subfigure}
	\begin{subfigure}[b]{0.5\textwidth}
		\centering
		\hspace{-25pt}\includegraphics[trim=3cm 3cm 3cm 0.1cm,scale=0.3,clip]{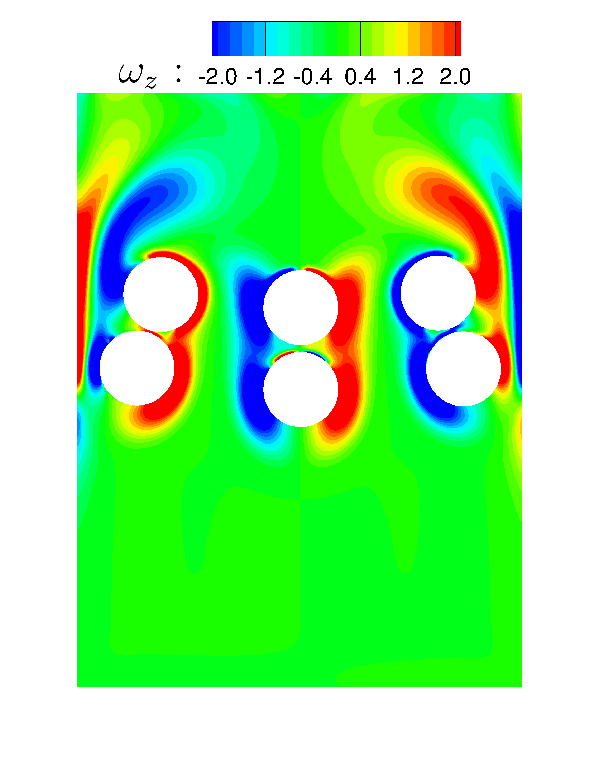}
		\caption{$t=9$}
		\label{fig:6fall3}
	\end{subfigure}%
	\begin{subfigure}[b]{0.5\textwidth}
		\centering
		\hspace{-25pt}\includegraphics[trim=3cm 3cm 3cm 0.1cm,scale=0.3,clip]{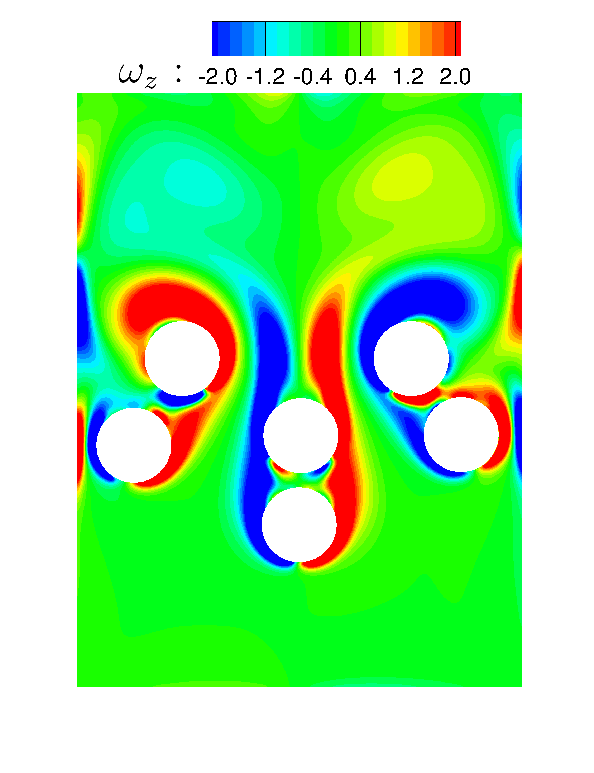}
		\caption{$t=12$}
		\label{fig:6fall4}
	\end{subfigure}
\caption{$\omega_z$ contour of six free falling particles (I)}
\label{fig:6falla}
\end{figure}
\begin{figure} \centering
	\begin{subfigure}[b]{0.5\textwidth}	
		\centering
		\hspace{-25pt}\includegraphics[trim=3cm 3cm 3cm 0.1cm,scale=0.3,clip]{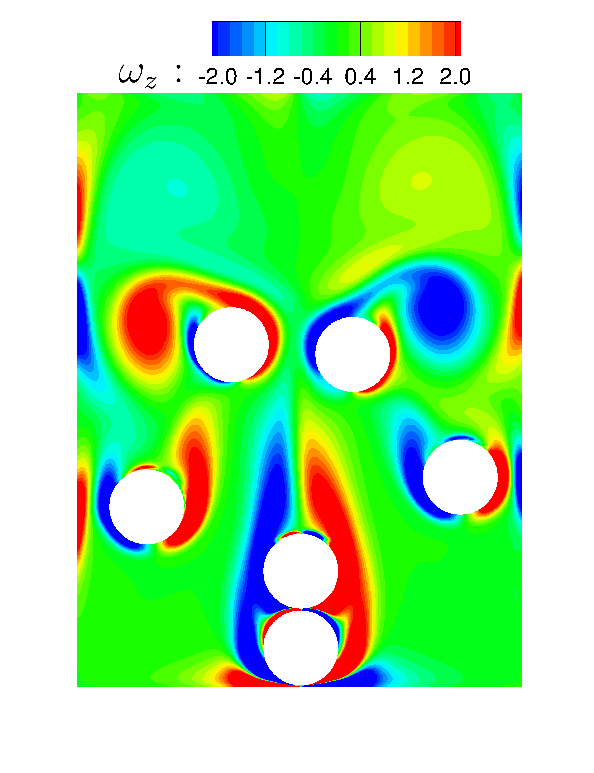}
		\caption{$t=15$}
		\label{fig:6fall5}
	\end{subfigure}%
	\begin{subfigure}[b]{0.5\textwidth}	
		\centering
		\hspace{-25pt}\includegraphics[trim=3cm 3cm 3cm 0.1cm,scale=0.3,clip]{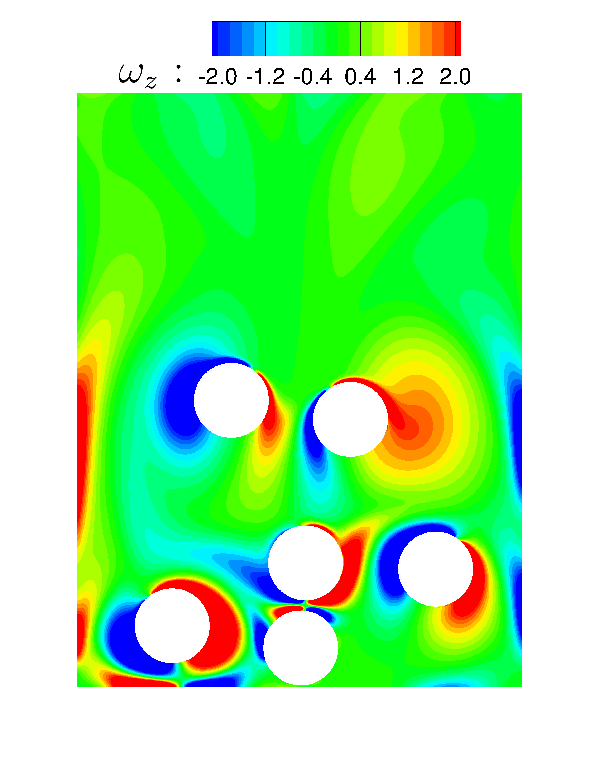}
		\caption{$t=18$}
		\label{fig:6fall6}
	\end{subfigure}
	\begin{subfigure}[b]{0.5\textwidth}
		\centering
		\hspace{-25pt}\includegraphics[trim=3cm 3cm 3cm 0.1cm,scale=0.3,clip]{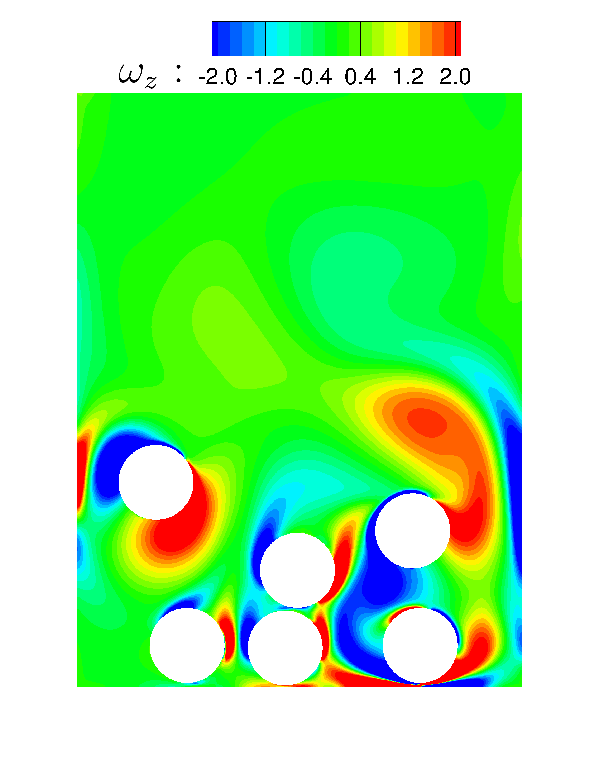}
		\caption{$t=21$}
		\label{fig:6fall7}
	\end{subfigure}%
	\begin{subfigure}[b]{0.5\textwidth}
		\centering
		\hspace{-25pt}\includegraphics[trim=3cm 3cm 3cm 0.1cm,scale=0.3,clip]{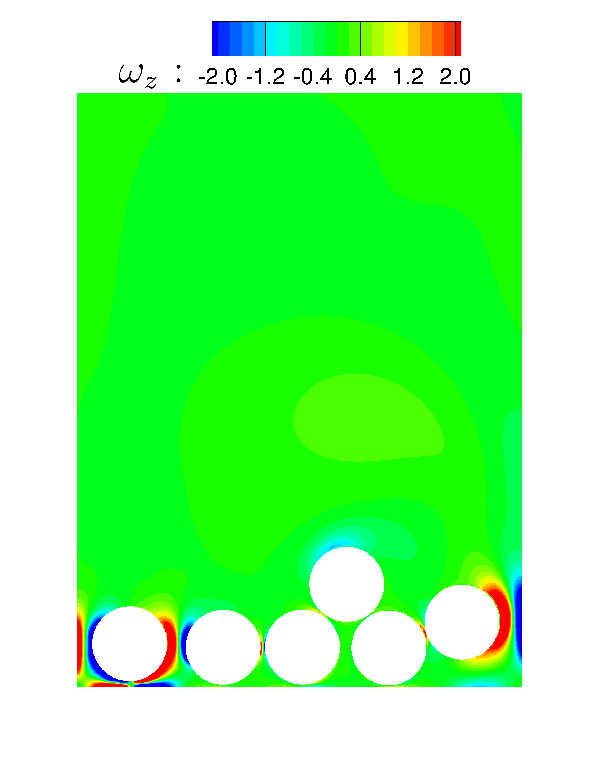}
		\caption{$t=24$}
		\label{fig:6fall8}
	\end{subfigure}
	\caption{$\omega_z$ contour of six free falling particles (II)}
	\label{fig:6fallb}
\end{figure}

In this benchmark example, six particles are freely falling under gravity from the rest. This problem is significantly differ from the single particle example in Sect.~\ref{sec:1fall}, because of the complex interaction between particles, walls and wakes. The objective is to demonstrate the robustness of proposed PGQ technique to handle much more challenging circumstances, e.g., rigid-body contact. Since we did not find literature of similar numerical or experimental setups, this example is meant to qualitatively demonstrate the capability of proposed PGQ technique. The width and height of domain are $x/D = [-3D,3D]$, $y/D=[1D,-7D]$ respectively, where $D=1.0$ is particle diameter. The top layer particles are rest at $x/D=0$ at $t=0$. The boundary conditions are identical to the benchmark example in Sect.~\ref{sec:1fall}. The fluid density, dynamic viscosity, mass ratio respectively are $\rho^f = 1.0$, $\mu=0.01$ and $m^* = 1.1$.

The implemented contact model~\cite{Wan2006IJfNMiF} ensures there is no penetration among particle and wall. Complex vortex wakes are generated as particle falling through the channel in Fig.~\ref{fig:6falla} and~\ref{fig:6fallb}. Eventually, all particles rest at the bottom of the channel and vortex wakes vanish.

\section{Conclusion} \label{sec:con}
A projection-based numerical integration technique, PGQ, was proposed for the application of FSI problems. This scheme is formulated based on tessellation technique. It operates on the matrix level, after the standard numerical integration rule, e.g., Gauss-Legendre Quadrature, is applied in each integration cell. 

Its main advantages are (1) no change in FE formulation and Quadrature rule for elements with/without embedded discontinuity, which simplifies implementation and improves its scalability to other physical problems, (2) variationally consistent with the derivation of FE formulation and well-suited for FE formulation, (3) approximation of the discontinuity with reduced dimension space. It possesses important characteristics: (1) \emph{partition of unity} property, (2) exact recovery of Gauss quadrature, (3) projection in quadratic form and (4) reduced-order modeling.

PGQ is implemented in various benchmark examples to assess its robustness and accuracy. It was shown the obtained numerical results via PGQ matched well with literature of various FSI applications. Therefore, the propose PGQ is excellent for numerical integration over cut cell with embedded discontinuities in FE framework for application of FSI problems.

\section*{Acknowledgments}
The first author would like to thank for the financial support from National Research Foundation through Keppel-NUS Corporate Laboratory. The conclusions put forward reflect the views of the authors alone, and not necessarily those of the institutions.

\bibliography{refs}

\end{document}